\documentclass[preprints,article,accept,moreauthors,pdftex]{my-mdpi}

\firstpage{1}
\pubvolume{11}
\issuenum{4}
\articlenumber{170}
\doinum{10.3390/axioms11040170}
\pubyear{2022}
\copyrightyear{2022}
\externaleditor{Academic Editor: Giampiero Palatucci}
\datereceived{31 October 2021}
\daterevised{25 March 2022} 
\dateaccepted{6 April 2022}
\datepublished{11 April 2022}
\hreflink{https://doi.org/10.3390/\linebreak axioms11040170} 

\pdfoutput=1

\newcommand{\alfa}{^{\alpha}}
\newcommand{\Dl}{{}^{\scriptscriptstyle C}_{\scriptscriptstyle 0}\!D^{\alpha}_{\scriptscriptstyle t}}

\Title{Fractional Modelling and Optimal Control of COVID-19 Transmission in Portugal}

\TitleCitation{Fractional Modelling and Optimal Control of COVID-19 Transmission in Portugal}

\Author{Silv\'{e}rio Rosa $^{1,\dagger}$\orcidA{}
and Delfim F. M. Torres $^{2,}$*$^{,\dagger}$\orcidB{}}

\AuthorNames{Silv\'{e}rio Rosa and Delfim F. M. Torres}

\AuthorCitation{Rosa, S.; Torres, D.F.M.}

\address{$^{1}$ \quad Department of Mathematics, Instituto de Telecomunicações (IT), Universidade da Beira Interior,
6201-001~Covilh\~{a}, Portugal; rosa@ubi.pt\\
$^{2}$ \quad Center for Research and Development in Mathematics and Applications (CIDMA),
Department of Mathematics, University of Aveiro, 3810-193 Aveiro, Portugal}

\corres{Correspondence: delfim@ua.pt}

\firstnote{These authors contributed equally to this work.}

\abstract{A fractional-order compartmental model was recently used
to describe real data of the first wave of the COVID-19 pandemic in Portugal
[Chaos Solitons Fractals 144 (2021), Art.~110652]. Here,
we modify that model in order to correct time dimensions
and use it to investigate the third wave of COVID-19
that occurred in Portugal from December 2020 to February 2021,
and that has surpassed all   previous waves, both in number and consequences.
A new fractional optimal control problem is then formulated and solved,
with vaccination and preventive measures as controls. A cost-effectiveness
analysis is carried out, and the obtained results   are discussed.}

\keyword{compartmental models; COVID-19 pandemic; third wave of COVID-19 in Portugal;
fractional-order calculus; optimal control}

\MSC{26A33; 34A08; 49N90; 92C60}

\begin{document}


\section{Introduction}

In January 2020, the World Health Organization (WHO) announced
the existence of a significant number of pneumonia cases in Wuhan.
Against all the predictions, COVID-19 (COrona VIrus Disease-19)
spread quickly across the globe and, on 11th of March,   was
declared as a pandemic \cite{WHO}. Caused by SARS-CoV-2
(Severe Acute Respiratory Syndrome Corona Virus 2), COVID-19
is the first pandemic in the digital era from which very few territories
of the world are untouched. Many governments were forced to decree
measures that seemed to be outdated, such as the isolation of individuals
and the complete lockdown of regions and even countries, that compromise
individual freedoms, damage business and economy, and threaten
a significant number of jobs.

To fight COVID-19 and its harmful effects,
a multidisciplinary approach is needed. In particular,
mathematical modelling plays an important role in  the prediction
of possible scenarios and in its effective control
{\cite{MyID:461,Bracher:et:al}.}
Readers interested in fractional modelling are referred to
\cite{Rev2:1,Rev2:2} and references therein.

The pandemic numbers have put national health systems
under pressure. Many reported cases were not reported on time,
but with a delay of days. Hence, in this paper, we do not use the number
of daily reported cases but the means of the previous five days
of reported cases, as suggested in \cite{NDAIROU2021110652}.

The mean of five days of daily reported cases induces memory into the model.
Fractional derivatives have been intensively used to obtain models of
infectious diseases that take into account the memory effects.
Many researchers have focused particular attention in  modelling
real-world phenomena using non-integer order derivatives. Those dynamics
have been modelled and studied by using the concept of fractional-order derivatives.
These problems appear in a range of diversified fields of applied sciences,
including biology, physics, ecology, and engineering \cite{MyID:476,MR1658022}.

A classical compartmental model with a super-spreader class was firstly
applied to give an estimation of infected subjects and fatalities in Wuhan, China,
in \cite{MR4093642}. The first-order derivative was then substituted
by a derivative of a fractional order, resulting in a model investigated
with Caputo fractional derivatives \cite{NDAIROU2021110652}. Here,
this fractional model is corrected and then used to model the third
wave of COVID-19 in Portugal. We start by the estimation of parameters
that best fit  the real data.  The sensitivity analysis to the fractional-order
model is performed in order to identify which model parameters are most
influential on the dynamics of the disease. Afterwards, fractional optimal
control is applied, showing the effectiveness of our approach.

This paper is organized as follows. In Section~\ref{2:section}, the fractional
order model is formulated. Our main results are then given in
Section~\ref{3:mainresults}: parameter estimation of the COVID-19 model
with real data of Portugal (Section~\ref{31:parameterestimation});
sensitivity analysis of the parameters of the model, without forgetting
the effect of the order of fractional differentiation
(Section~\ref{subsec:sensitiveanalysis}); fractional optimal control
of the model (Section~\ref{subsec:control}); and, finally, numerical
simulations and cost-effectiveness of the fractional model
(Section~\ref{subsec:numres}). We end with Section~\ref{sec:conc},
which states the conclusions. Our main theoretical contributions
consist of fractional-order model consistency and the mathematical
problem rearrangement according to the Pontryagin theory. Moreover,
we solve the fractional optimal control problem numerically
by a method of Adams--Basforth--Moulton.


\section{Fractional-Order COVID-19 Model}
\label{2:section}

Since the time interval considered in this study is small, we assume
that population is invariant. In addition, population is divided in eight
classes: (i) susceptible individuals, \emph{S}; (ii) exposed individuals, \emph{E};
(iii) symptomatic infectious individuals, \emph{I}; (iv) super-spreader  individuals, \emph{P};
(v) asymptomatic infectious individuals, \emph{A}; (vi) hospitalized individuals, \emph{H};
(vii) recovered individuals, \emph{R}; and (viii) fatality class, \emph{F}.
Our fractional-order model is derived from the one presented in
Nda\"{\i}rou et al. \cite{NDAIROU2021110652}, which gives a generalization
of an integer-order model that was used to study the start of the pandemic
in  Wuhan \cite{MR4093642}. We use the fractional derivative
in the sense of Caputo. Fractional differential equations are an active area
of research and are adequate to incorporate the history of the processes.
The model system of equations for COVID-19 proposed in \cite{NDAIROU2021110652}
is given by
\begin{equation}
\label{Covid_model}
\begin{cases}
\Dl S = & -\beta\dfrac{I S}{N}-l \beta
\dfrac{H S}{N} -\beta' \dfrac{P S}{N},\\[2mm]
\Dl E = & \beta\dfrac{I S}{N}+l \beta
\dfrac{H S}{N} +\beta' \dfrac{P S}{N}-\kappa E,\\[2mm]
\Dl I = & \kappa \rho_1 E -(\gamma_a+\gamma_i)I-\delta_i I,\\[1mm]
\Dl P = & \kappa \rho_2 E -(\gamma_a+\gamma_i) P -\delta_p P,\\[1mm]
\Dl A = & \kappa(1-\rho_1-\rho_2)E,\\[1mm]
\Dl H = & \gamma_a(I+P)-\gamma_r H-\delta_h H,\\[1mm]
\Dl R = & \gamma_i(I+P)+\gamma_r H,\\[1mm]
\Dl F = & \delta_i I+\delta_p P +\delta_h H,
\end{cases}
\end{equation}
where $\Dl$ denotes the left Caputo fractional-order derivative
with derivative order  $\alpha$ ($0<\alpha \leqslant 1$).
A description of the parameters of Model \eqref{Covid_model}
can be found in Table~\ref{tab:param}.
\begin{table}[H]
\caption{Description and values of the parameters
of Model \eqref{Covid_model} taken from \cite{MR4093642,MyID:460}.}
\setlength{\tabcolsep}{5.8mm}{
\begin{tabular}{clc}
\toprule
\textbf{Name} & \multicolumn{1}{c}{\textbf{Description}} & \textbf{Value} \\
\midrule
$\beta$ &\mbox{human-to-human transmission coefficient}&2.55 \\
$l$ &\text{transmissibility of hospitalized patients}& 1.56 \\
$\beta'$ &  \text{transmission coefficient of super-spreaders}& 7.65\\
\multirow{2}{*}{$\kappa$}
&\text{rate at which an individual leaves the exposed}&	\multirow{2}{*}{ 0.25}\\
&\text{class to become infectious}& \\
\multirow{2}{*}{$\rho_1$}
&\text{proportion of progression from class $E$ }& \multirow{2}{*}{0.58}  \\
&\text{to symptomatic infectious class $I$}& \\
$\rho_2$ &\text{rate at which exposed ind. become super-spreaders} & 0.001  \\
\multirow{2}{*}{$\gamma_a$} & \text{rate at which symptomatic and super-spreaders}
& \multirow{2}{*}{0.94}  \\
&\text{ become hospitalized }&\\
$\gamma_i$ & \text{recovery rate without being hospitalized}& 0.27  \\
$\gamma_r$ & \text{recovery rate of hospitalized patients}& 0.5  \\
$\delta_i$ & \text{disease induced death rate due to infected ind.}& 1/23  \\
$\delta_p$ & \text{disease induced death rate due to super-spreader  ind.}& 1/23  \\
$\delta_h$ & \text{disease induced death rate due to hospitalized  ind.}& 1/23  \\
\bottomrule
\end{tabular}}
\label{tab:param}
\end{table}

We note that the equations of Model \eqref{Covid_model}
do not have appropriate time dimensions.
Indeed, on the left-hand side, dimension is (time)$^{-\alpha}$, while
on the right-hand side, dimension is (time)$^{-1}$.
This means that Model \eqref{Covid_model} is only consistent
when $\alpha = 1$. For the importance to be consistent with dimensions,
we refer the reader, e.g., to \cite{MR3808497,MR3928263}. Thus,
here, we correct System \eqref{Covid_model} as follows:
\begin{equation}
\label{Covid_model:exp}
\begin{cases}
\Dl S = & -\beta\alfa\dfrac{I S}{N}-l \beta\alfa
\dfrac{H S}{N} -{\beta'}\alfa \dfrac{P S}{N},\\[2mm]
\Dl E = & \beta\alfa\dfrac{I S}{N}+l \beta\alfa
\dfrac{H S}{N} +{\beta'}\alfa \dfrac{P S}{N}-\kappa\alfa E,\\[2mm]
\Dl I = & \kappa\alfa \rho_1 E -(\gamma_a\alfa+\gamma_i\alfa)I
-\delta_i\alfa I,\\[1mm]
\Dl P = & \kappa\alfa \rho_2 E -(\gamma_a\alfa+\gamma_i\alfa) P
-\delta_p\alfa P,\\[1mm]
\Dl A = & \kappa\alfa(1-\rho_1-\rho_2)E,\\[1mm]
\Dl H = & \gamma_a\alfa(I+P)-\gamma_r\alfa H-\delta_h\alfa H,\\[1mm]
\Dl R = & \gamma_i\alfa(I+P)+\gamma_r\alfa H,\\[1mm]
\Dl F = & \delta_i\alfa I+\delta_p\alfa P +\delta_h\alfa H.
\end{cases}
\end{equation}

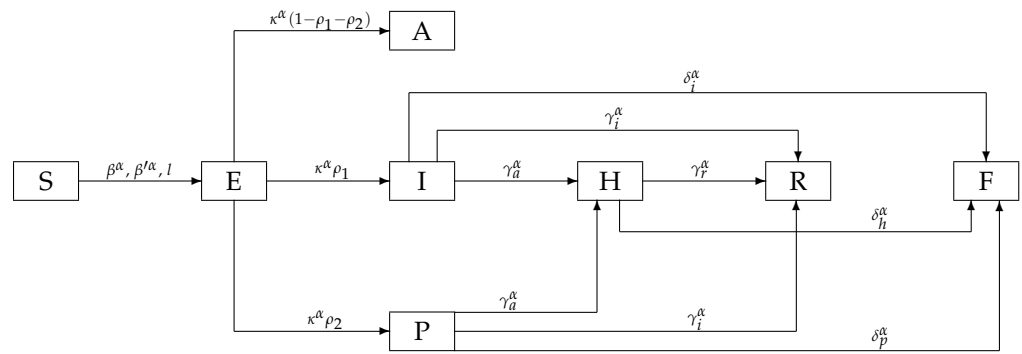
\begin{figure}[H]
\begin{picture}(560,145)(0,-20)
\setlength{\unitlength}{.25mm}
\put(0,60){\framebox(34,20){S}}
\put(200,140){\framebox(34,20){A}}
\put(100,60){\framebox(34,20){E}}
\put(200,-20){\framebox(34,20){P}}
\put(200,60){\framebox(34,20){I}}
\put(300,60){\framebox(34,20){H}}
\put(400,60){\framebox(34,20){R}}
\put(500,60){\framebox(34,20){F}}
\put(34,70){\vector(1,0){66}}
\put(49,73){$\scriptscriptstyle \beta^{\alpha},\;\beta'^{\alpha},\; l$}
\put(134,70){\vector(1,0){66}}
\put(160,74){$\scriptscriptstyle \kappa^{\alpha}\rho_1$}
\put(234,70){\vector(1,0){66}}
\put(260,74){$\scriptscriptstyle \gamma_a^{\alpha}$}
\put(334,70){\vector(1,0){66}}
\put(360,74){$\scriptscriptstyle \gamma_r^{\alpha}$}
\put(117,80){\line(0,1){70}}
\put(117,150){\vector(1,0){83}}
\put(136,153.5){$\scriptscriptstyle \kappa^{\alpha}(1-\rho_1-\rho_2)$}
\put(210,80){\line(0,1){37}}
\put(210,117){\line(1,0){307}}
\put(517,117){\vector(0,-1){37}}
\put(356,122){$\scriptscriptstyle \delta_i^{\alpha}$}
\put(225,80){\line(0,1){17}}
\put(225,97){\line(1,0){192}}
\put(417,97){\vector(0,-1){17}}
\put(314,102){$\scriptscriptstyle \gamma_i^{\alpha}$}
\put(117,60){\line(0,-1){70}}
\put(117,-10){\vector(1,0){83}}
\put(156,-6){$\scriptscriptstyle \kappa^{\alpha}\rho_2$}
\put(322,60.1){\line(0,-1){17}}
\put(322,43){\line(1,0){187}}
\put(509,43){\vector(0,1){17}}
\put(456,49){$\scriptscriptstyle \delta_h^{\alpha}$}
\put(234,0.1){\line(1,0){76}}
\put(310,0.1){\vector(0,1){60}}
\put(258,4){$\scriptscriptstyle \gamma_a^{\alpha}$}
\put(234,-10){\line(1,0){182}}
\put(416,-10){\vector(0,1){70}}
\put(358,-5){$\scriptscriptstyle \gamma_i^{\alpha}$}
\put(234,-20.3){\line(1,0){290}}
\put(524,-20.3){\vector(0,1){80}}
\put(456,-15){$\scriptscriptstyle \delta_p^{\alpha}$}
\end{picture}
\caption{Flow diagram of the disease dynamics according
to Model \eqref{Covid_model:exp}.}\label{diagram:model}
\end{figure}

A flow diagram of Model \eqref{Covid_model:exp}
is given in Figure~\ref{diagram:model}.
Note that, in the particular case $\alpha = 1$, both
Models \eqref{Covid_model} and \eqref{Covid_model:exp}
coincide with the classical COVID-19 model
first introduced in~\cite{MR4093642}.
System \eqref{Covid_model:exp} will be investigated
in Section~\ref{3:mainresults}.
\vspace{-6pt}


\section{Main Results}
\label{3:mainresults}

We begin by discussing the adherence of the corrected
Model \eqref{Covid_model:exp} introduced
in \mbox{Section~\ref{2:section}} with respect to COVID-19
and real data from the third wave that occurred in Portugal.
For that, we need an adequate estimation of the parameter values.


\subsection{Parameter Estimation}
\label{31:parameterestimation}

The uncorrected model \eqref{Covid_model} was used to study the dynamics of COVID-19
at the early stage of  the pandemic, between 3 March 2020 and 27 April 2020,
in \cite{NDAIROU2021110652}. During that period, government decreed
a general lockdown of the country that was well accepted by the population
due to the impact of the disease in other countries where it started earlier.
With limited knowledge, this revealed to be an effective measure to reduce
contacts and control the disease in a relatively short period of time.

That initial stage ended in May 2020, with the gradual release of the country
from COVID-19 container measures and  the cancelling of the State of Emergency.
In that date, preventive measures were adopted to control the disease.
One of the most important was the use of masks in confined spaces,
made mandatory by Portuguese \emph{Decreto-Lei n.º 20/2020}.

In October 2020, due to autumn weather conditions, the number of infected
individuals started to be worrisome. In order to control it and to protect
the Portuguese National Health Service, in November 2020, a new series
of States of Emergency began.  In January 2021, the cold weather,
some relaxation among  the population concerning preventive measures
(being in crowds or poorly ventilated spaces,   the misuse of masks, \ldots)
combined with new variants of COVID-19---more contagious than ever---
that started circulating, forced the government to declare a new lockdown
with the closure of schools in 22 January 2021.

\begin{figure}[H]
\includegraphics[scale=0.7]{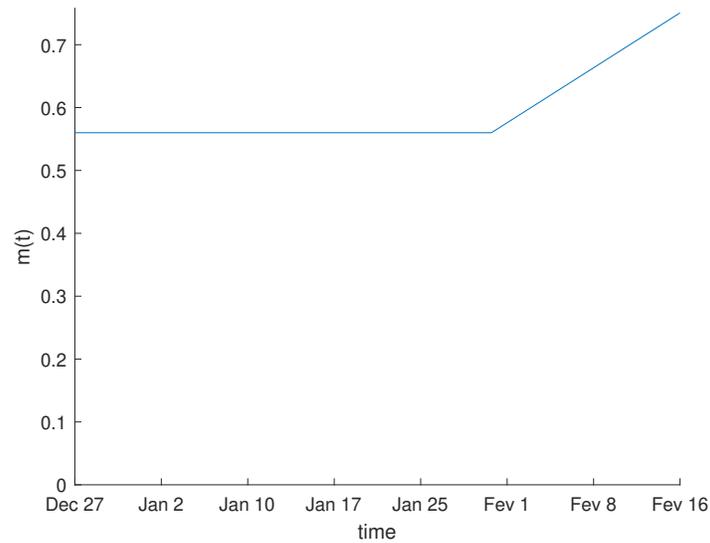}
\caption{Evolution function $m(t)$.}
\label{fig:function_m}
\end{figure}

Due to the implementation of the preventive measures, during our model fitting,
we assume that parameters have the same values of the first wave, with  the exception
of contact rates. The period we chose starts in 27 December 2020 and ends
16 February 2021, covering the third wave of the pandemic in Portugal.
During that period, the closure of schools was declared, and that most certainly had
an impact on contact rates. Hence, we consider that
the transmission coefficient $\beta$ is replaced by $\beta(1- m(t))$,
and that $\beta'$ is replaced by $\beta'(1- m(t))$, where $m(t)$ is a
continuous function that represents the rate of  reduction of the contacts,
and that varies with time according to Figure~\ref{fig:function_m}.
\vspace{-6pt}

The values of two parameters were  determined by the fitting of the model:
(i) derivative order of the model, $\alpha$,  and (ii) a scaling factor $s$.
See Table~\ref{tab:fitting} for the resulting values and errors
associated with the fitting.
\begin{table}[H]
\caption{Results of model fitting.}
\setlength{\tabcolsep}{6.4mm}{
\begin{tabular}{cccc}
\toprule
\multicolumn{1}{c}{\textbf{Derivative Order}}  &\boldmath{$s$}
&  \textbf{Absolute Error}  & \textbf{Relative Error (\%)}\\
\midrule
1.0  & 21.08&8595& 14.13\\
0.99  & 19.87& 8135&13.37\\
\bottomrule
\end{tabular}}
\label{tab:fitting}
\end{table}

\begin{figure}[H]
\includegraphics[scale=0.5]{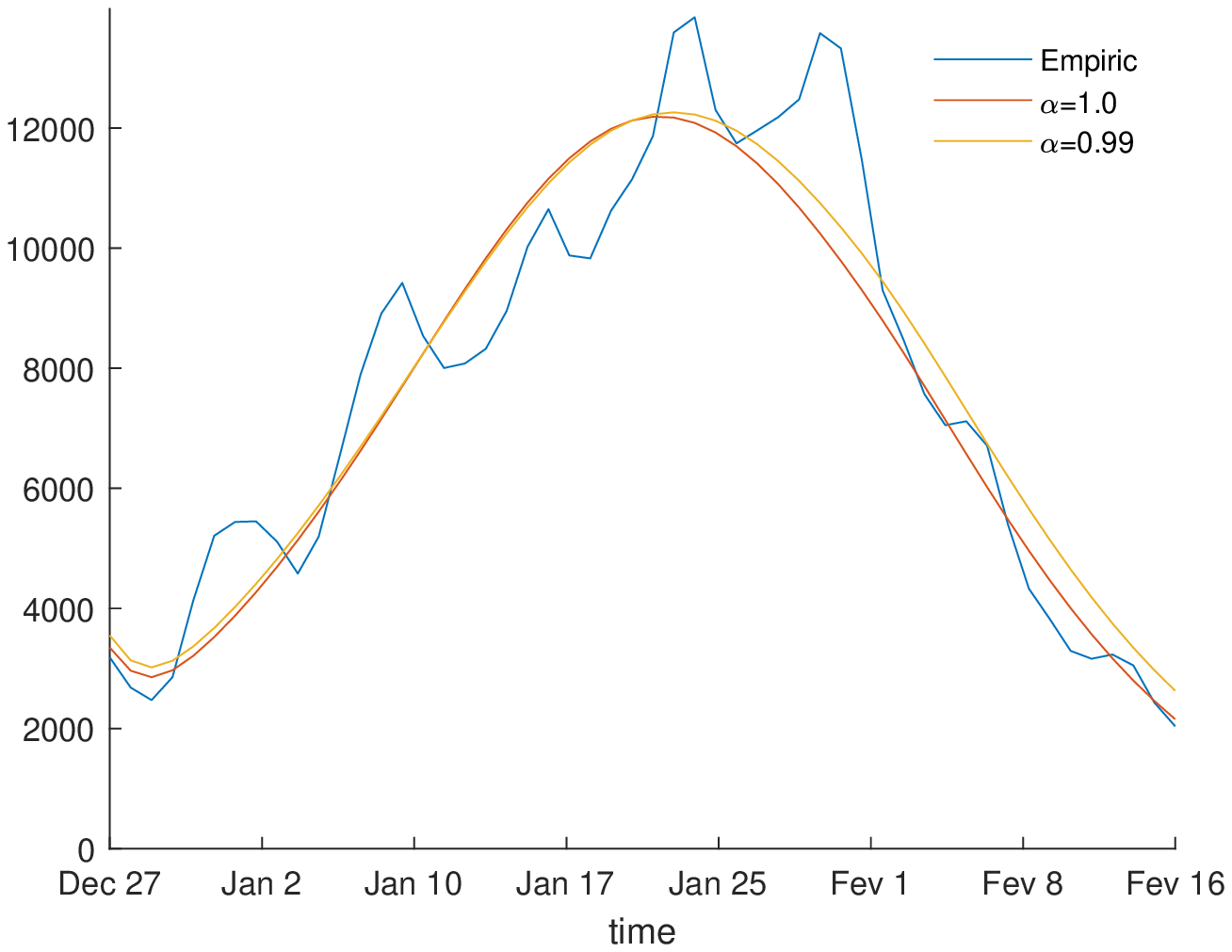}
\includegraphics[scale=0.5]{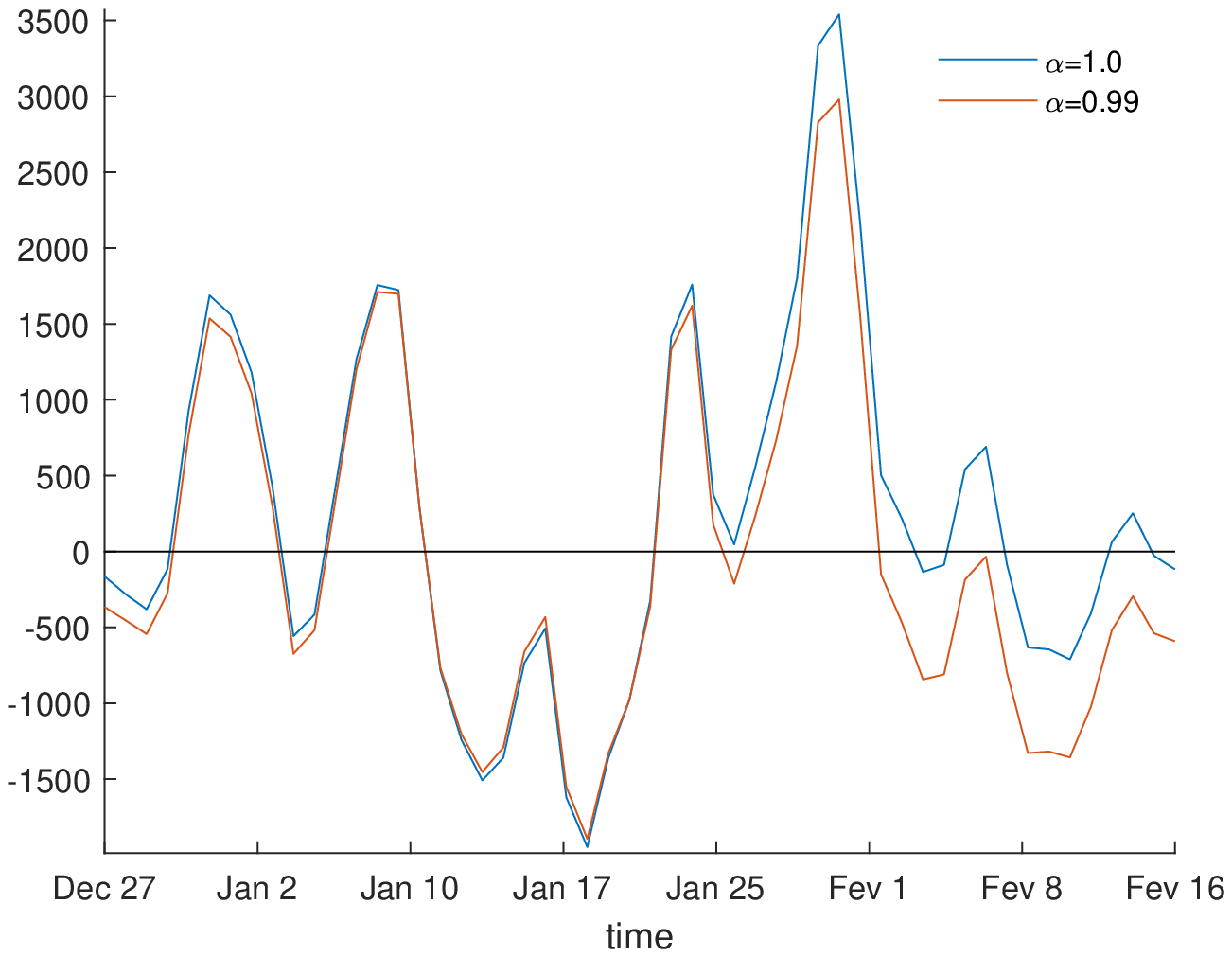}
\caption{(\textbf{left}) The number of confirmed cases per day in Portugal
versus the ones predicted by Model \eqref{Covid_model:exp}
with parameters given by Table~\ref{tab:param}. The blue line
corresponds to the real data ($I+P+H$) and the remaining lines
have been obtained solving numerically the system of fractional
differential equations \eqref{Covid_model:exp}. (\textbf{right})
The difference between the number of confirmed cases per day and the number
of estimated cases,  the solution of \eqref{Covid_model:exp}.}
\label{fig:gen_infected_predicted_estimated}
\end{figure}

The fitting curves are presented
in Figure~\ref{fig:gen_infected_predicted_estimated},
where real data was obtained from \cite{DGS21,Data}.

It is known that the available data has some mistakes. Frequently,
days reporting few cases are followed by days reporting many
new cases, without correspondence with reality. Due to the stress
that the pandemic imposed over doctors and health professionals,
these data do not flow as quickly and effectively as expected.
Therefore, the numbers of new cases are not correctly determined by
reported daily numbers. Thus, the number of new cases considered
in this manuscript is, as suggested in \cite{NDAIROU2021110652},
the mean of the previous five days of reported cases.

Following \cite{MR3771538}, our data fitting consisted
in minimization of the $l_2$ norm of the difference between
the real values and predictive cases of COVID-19 infection
given by \mbox{Model \eqref{Covid_model:exp}}. Due to the oscillation
of the number of confirmed cases per day, the fitting curves
observed in Figure~\ref{fig:gen_infected_predicted_estimated}
(left) have significant gaps in comparison with real data.
Figure~\ref{fig:gen_infected_predicted_estimated} (right) presents
the difference between the number of confirmed cases per day and the number
of estimated cases, showing that the proposed model approximates well
the average number of cases.

Consequently, fitting errors are quite high, being 14.13\%
for the classical integer-order model ($\alpha = 1$)
and 13.37\% for the fractional-order model with $\alpha=0.99$,
according to Table~\ref{tab:fitting}. The difference between
the two derivative orders, in terms of absolute
errors, justifies the preference for the fractional-order model
with respect to the classical one.


\subsection{Sensitivity Analysis}
\label{subsec:sensitiveanalysis}

An important threshold, while studying infectious disease models,
is the basic reproduction number. Following \cite{MR4093642,MyID:460},
we conclude that the basic reproduction number of the COVID-19 model
\eqref{Covid_model:exp} is
\begin{equation}
\label{r0_model}
{R_0}=\frac{\beta\alfa \rho_1(\gamma_a\alfa l +a_h)}{a_i a_h}
+\frac{(\beta\alfa \gamma_a\alfa l +{\beta'}\alfa a_h)\rho_2}{a_p a_h},
\end{equation}
where $a_i= \gamma_a\alfa+\gamma_i\alfa+\delta_i\alfa$,
$a_p= \gamma_a\alfa+\gamma_i\alfa+\delta_p\alfa$
and $a_h= \gamma_r\alfa+\delta_h\alfa$.

\begin{figure}[H]
\includegraphics[scale=0.55]{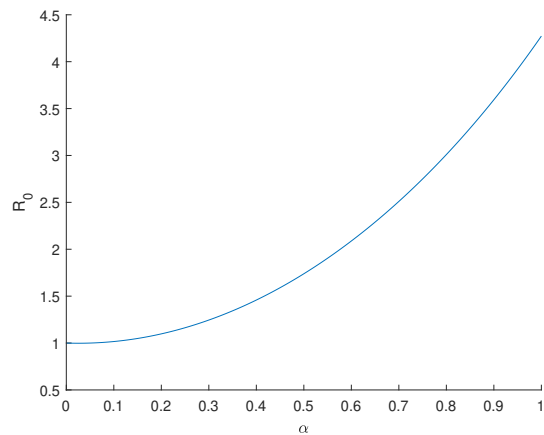}
\caption{Impact of the variation of the derivative order, $\alpha$,
in the evolution of the basic reproduction number ${R_0}$
of the COVID-19 model \eqref{Covid_model:exp}.}
\label{fig:R0_varalpha}
\end{figure}

The impact of the variation of the derivative order $\alpha$,
in the evolution of ${R_0}$, is presented in
Figure~\ref{fig:R0_varalpha}. We observe in Figure~\ref{fig:R0_varalpha}
that $R_0\geqslant 1$, that is, we have an endemic scenario
regardless of the value of $\alpha$.

Sensitivity analysis measures the importance of each parameter
of the model in the disease transmission.

\begin{Definition}[The sensitivity index
\cite{chitnis2008determining,rodrigues2013sensitivity}]
\label{def:sentInd}
Let ${R_0}$ be differentiable with respect to a given parameter $p$.
The sensitivity of the model with respect to that parameter $p$
can be measured by the index $\Upsilon_p^{{R_0}}$ given by
$$
\Upsilon_p^{{R_0}}=\frac{\partial {R_0}}{\partial p}\frac{p}{{R_0}}.
$$
\end{Definition}

The sensitivity analysis of the classical COVID-19 model
was presented in \cite{MR4093642}, that is, in the
particular situation $\alpha = 1$ in \eqref{Covid_model:exp}.
In it, the authors concluded that the most sensitive parameters
to the basic reproduction number $R_0$ are $\beta$, $\rho_1$, and $l$.
Therefore, special attention should be paid to the estimation of those parameters.
In contrast, the estimation of $\beta'$, $\rho_2$, $\gamma_a$, $\delta_i$, $\delta_p$,
and $\delta_h$ does not require as much attention because of its low sensitivity.

In the general situation of our fractional-order model \eqref{Covid_model:exp},
it should be emphasized that the sensitivity depends on the derivative order $\alpha$
of the fractional operator. We can observe this in Figure~\ref{fig:sensitivity}, where
one clearly sees that the variation of $\alpha$ influences the sensitivity
index of parameters $\beta$, $\rho_1$, $\gamma_i$, $\gamma_r$, and $l$.
\begin{figure}[H]
\includegraphics[scale=0.5]{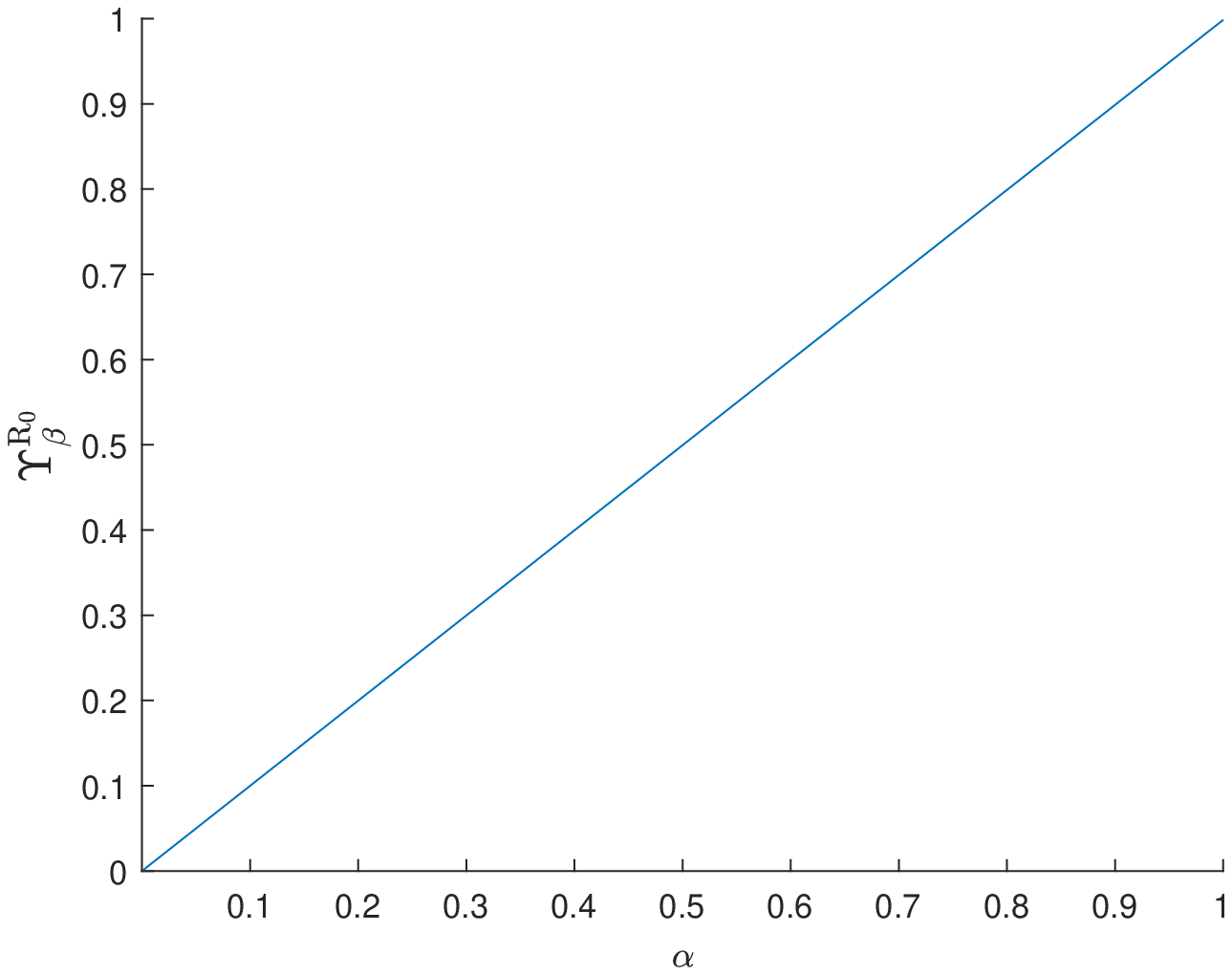}
\includegraphics[scale=0.5]{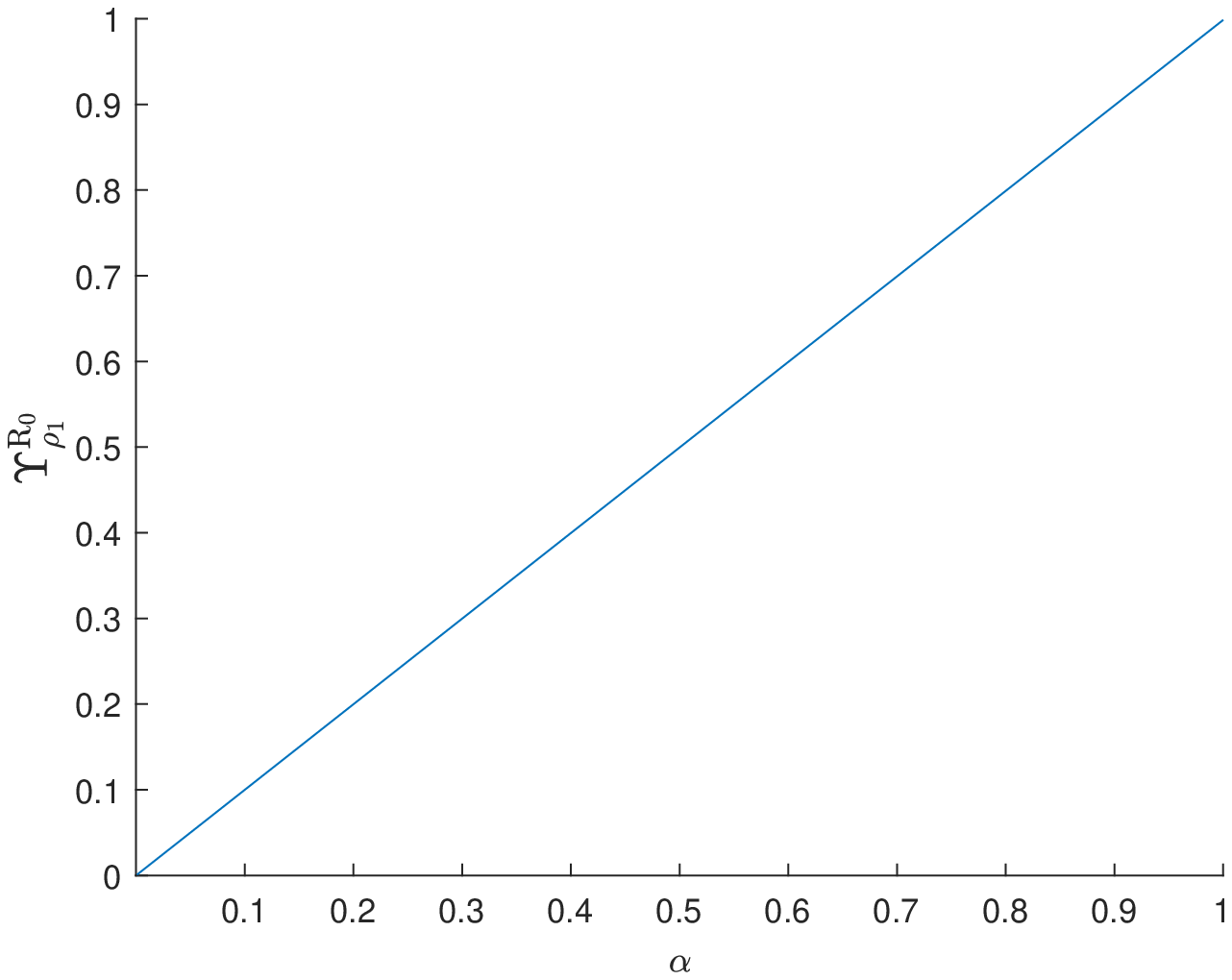}
\caption{\textit{Cont}.}
\label{fig:sensitivity}
\end{figure}

\begin{figure}[H]\ContinuedFloat
\includegraphics[scale=0.5]{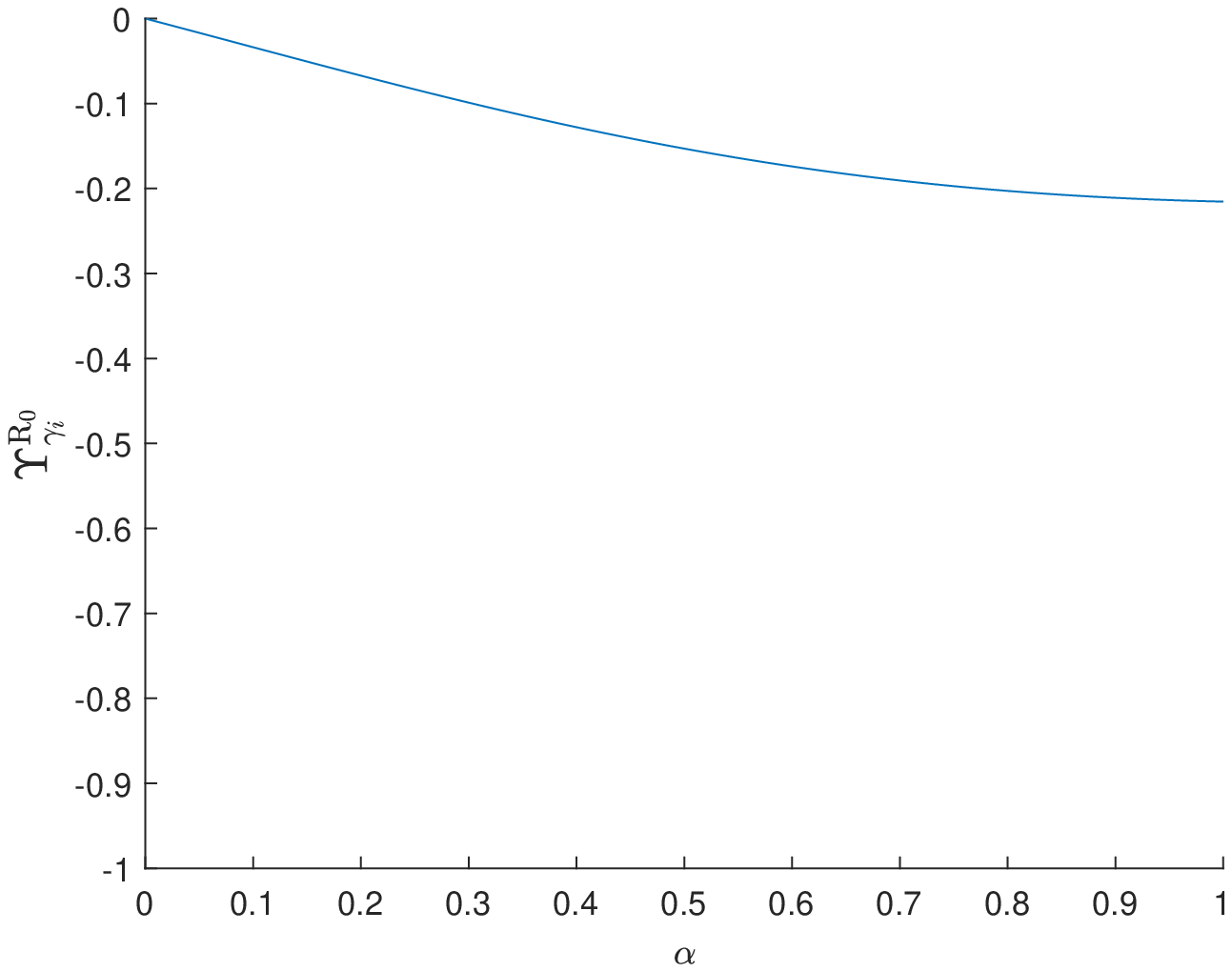}
\includegraphics[scale=0.5]{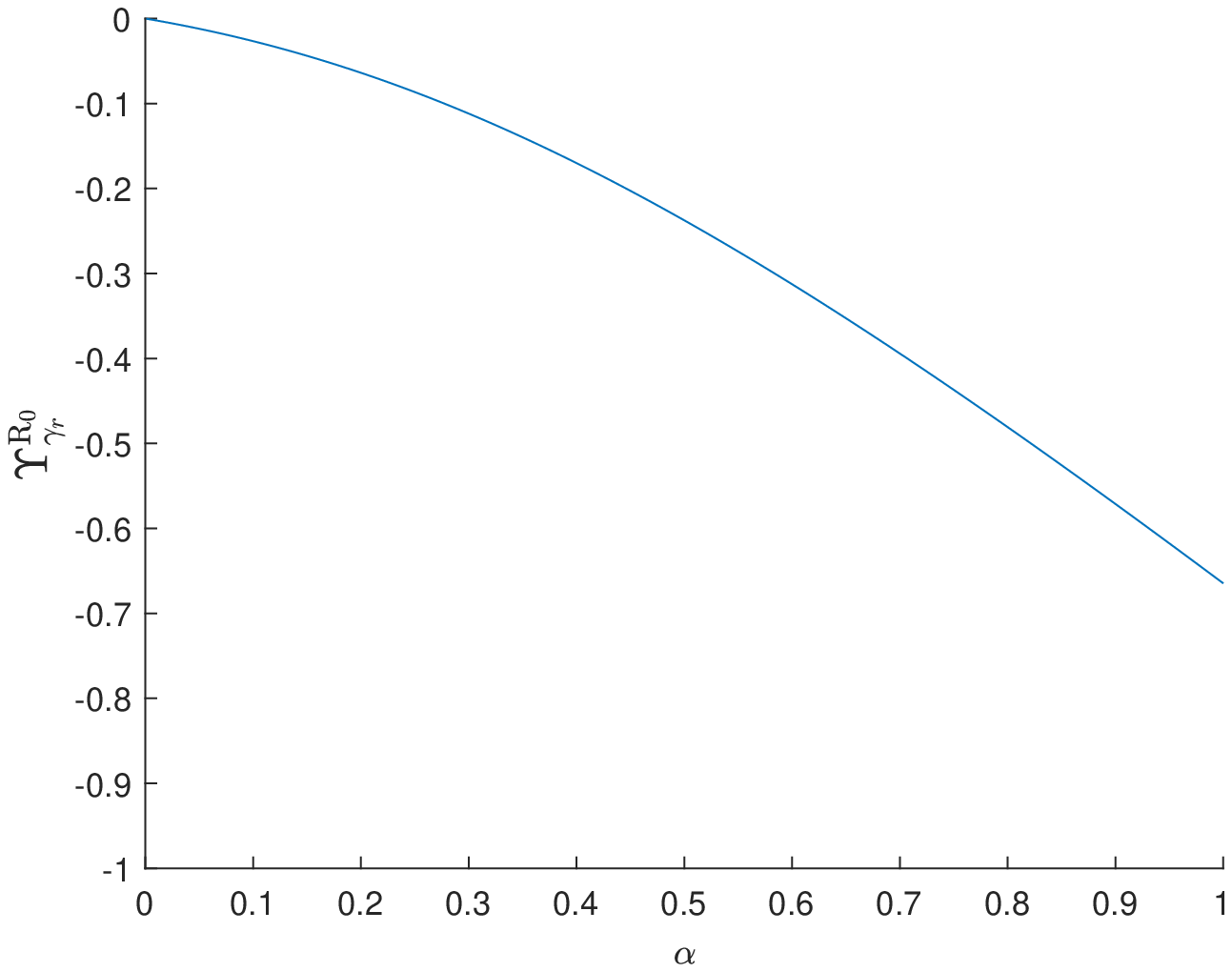}

\includegraphics[scale=0.5]{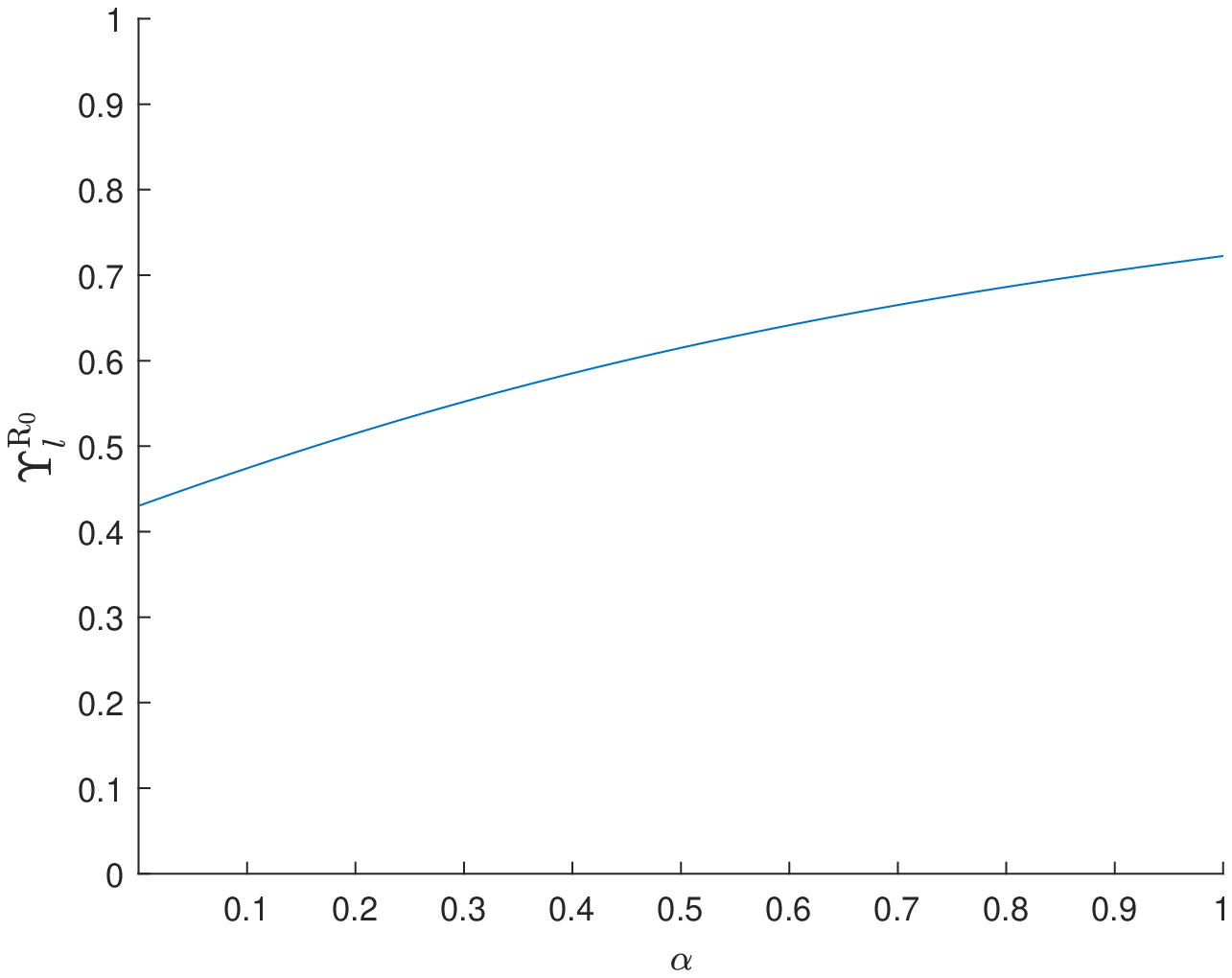}
\caption{Impact of the variation of $\alpha$ in the sensitivity
indexes of $\beta$, $\rho_1$, $\gamma_i$, $\gamma_r$, and $l$,
in agreement with Definition~\ref{def:sentInd}.}
\label{fig:sensitivity}
\end{figure}

We see in Figure~\ref{fig:sensitivity} that the sensitivity indexes
of $\beta$ and $\rho_1$ exhibit a quite similar evolution with respect
to $\alpha$, being very sensitive to the variation of $\alpha$. On the
opposite side, the sensitivity indexes of $\gamma_i$ and $l$ are much
less sensitive to the variation of the fractional-order $\alpha$.
The graphics with respect to the remaining parameters of the model
are omitted here because, using the same scale, their curves do not go
far from the $x$-axis.

\begin{figure}[H]
\includegraphics[scale=0.55]{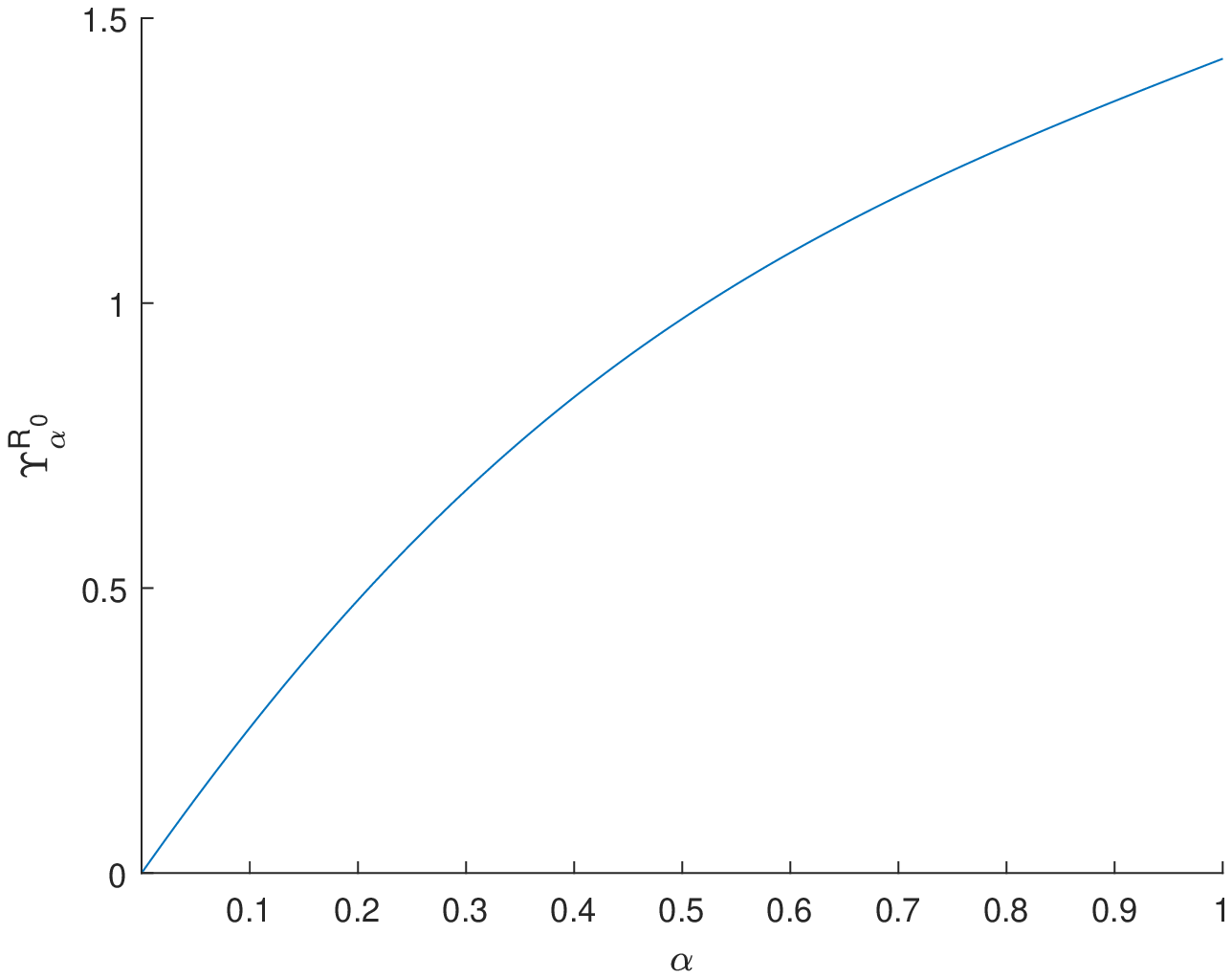}
\caption{The sensitivity of Model \eqref{Covid_model:exp}
with respect to the fractional order of differentiation $\alpha$,
in agreement with Definition~\ref{def:sentInd}.}
\label{fig:sensitivity_BRN}
\end{figure}

The evolution of the sensitivity index for the basic reproduction number $R{_0}$
with the variation of the derivative order $\alpha$ is presented
in Figure~\ref{fig:sensitivity_BRN}. We observe that the sensitivity index
decreases with the decrease of $\alpha$.

Figures~\ref{fig:sensitivity} and \ref{fig:sensitivity_BRN} also
show that the sensitivity of the model decreases in absolute value
with the decrease in the derivative order. This means that the fractional-order
model is less sensitive than the classical one, which is an advantage of our model.


\subsection{Fractional Optimal Control of the Model}
\label{subsec:control}

In this section, we aim to minimize the number of COVID-19-infected
individuals and reduce the cost of control measures. This is
carried out through  (i) the use of vaccination as an effective
measure time-dependent control $v(t)$ and (ii) the use of the preventive
measure control $m(t)$ (representing the use of masks, limitations
of the number of individuals in closed confined spaces, etc.),
in order to control person-to-person contact. Therefore, the following
fractional-order optimal control problem is considered:
\begin{align}
\label{cost-functional}
\begin{split}
\min ~\mathcal{J}(I, P,v,&m)
=\int_0^{t_f} \left(k_1 I+ k_2 P+k_3 v^2+k_4 m^2\right) ~dt\\
\end{split}\end{align}
subject to
\begin{equation}
\label{sys:state}
\begin{cases}
\Dl S =  -\beta\alfa(1-m)\dfrac{I S}{N}-l \beta\alfa(1-m) \dfrac{H S}{N}
-{\beta'}\alfa(1-m) \dfrac{P S}{N}-v S,\\[2mm]
\Dl E =  \beta\alfa(1-m)\dfrac{I S}{N}+l \beta\alfa(1 -m)\dfrac{H S}{N}
+{\beta'}\alfa(1 -m)\dfrac{P S}{N}-\kappa\alfa E,\\[2mm]
\Dl I =  \kappa\alfa \rho_1 E -(\gamma_a\alfa+\gamma_i\alfa)I-\delta_i\alfa I,\\[1mm]
\Dl P =  \kappa\alfa \rho_2 E -(\gamma_a\alfa+\gamma_i\alfa) P -\delta_p\alfa P,\\[1mm]
\Dl A =  \kappa\alfa(1-\rho_1-\rho_2)E,\\[1mm]
\Dl H =  \gamma_a\alfa(I+P)-\gamma_r\alfa H-\delta_h\alfa H,\\[1mm]
\Dl R =  \gamma_i\alfa(I+P)+\gamma_r\alfa H+ v S,\\[1mm]
\Dl F =  \delta_i\alfa I+\delta_p\alfa P +\delta_h\alfa H,
\end{cases}
\end{equation}
with given initial conditions
\begin{equation}
\label{ocp:ic}
\left(S(0),E(0),I(0),P(0),A(0),H(0),R(0),F(0)\right)
= \left(S_0,E_0,I_0,P_0,A_0,H_0,R_0,F_0\right) \geqslant 0.
\end{equation}

Parameters $0<k_1,\ldots,k_4 <+\infty$ are positive numbers that balance
the size of the linear and quadratic terms in the cost, and $t_f$ is the duration
of the control program. Moreover, $k_3$ and $k_4$ represent the costs
of applying the control measures $v$ and $m$, respectively.
The set of admissible control functions is
\begin{equation}
\label{omega:set}
\mathcal{U}=\left\lbrace (v(\cdot),m(\cdot))
\in L^{\infty}(0,t_{f}) : 0 \leqslant v(t)
\leqslant {v}_{\max},
\  0 \leqslant m(t) \leqslant {m}_{\max} \ \forall \ t \in (0,t_{f})
\right\rbrace.
\end{equation}

Pontryagin's minimum principle (PMP) for fractional optimal control
\cite{MR3443073} is used to determine the solution of the problem.
The Hamiltonian associated with our optimal control problem is given by
\begin{equation}
\begin{split}
H&= k_1 I+ k_2 P+k_3 v^2+k_4 m^2+\xi_1\bigg( -\beta\alfa(1-m)
\dfrac{I S}{N}-l \beta\alfa(1-m) \dfrac{H S}{N}\\
&-{\beta'}\alfa(1-m) \dfrac{P S}{N}-v S\bigg)+\xi_2\bigg(\beta\alfa(1-m)
\dfrac{I S}{N}+l \beta\alfa(1 -m)\dfrac{H S}{N}\\
&+{\beta'}\alfa(1 -m)\dfrac{P S}{N}-\kappa\alfa E\bigg)
+\xi_3\bigg(\kappa\alfa \rho_1 E -(\gamma_a\alfa+\gamma_i\alfa)I-\delta_i\alfa I\bigg)\\
&+\xi_4\bigg(\kappa\alfa \rho_2 E -(\gamma_a\alfa+\gamma_i\alfa) P -\delta_p\alfa
P\bigg)+\xi_5\bigg(\kappa\alfa(1-\rho_1-\rho_2)E\bigg)\\
&+\xi_6\bigg(\gamma_a\alfa(I+P)-\gamma_r\alfa H-\delta_h\alfa H\bigg)
+\xi_7\bigg(\gamma_i\alfa(I+P)+\gamma_r\alfa H+ v S\bigg)\\
&+\xi_8\bigg(\delta_i\alfa I+\delta_p\alfa P +\delta_h\alfa H\bigg)
\end{split}
\end{equation}
and the adjoint system of PMP asserts that the co-state variables
$\xi_i(t)$, $i = 1, \ldots, 8$, satisfy
\begin{equation}
\label{eq:co_states_fr2}
\begin{cases}
\Dl  \xi_1(t') = \frac{(m-1)(\beta\alfa(I+l H)+{\beta'}\alfa P)(\xi_1-\xi_2)}{N}
+(\xi_7-\xi_1)v, \\[1.2mm]
\Dl \xi_2(t') = \kappa\alfa(-\xi_2+\xi_3\rho_1
+\xi_4\rho_2-\xi_5(\rho_1+\rho_2-1)),\\[1.2mm]
\Dl \xi_3(t') = k_1-(\gamma_a\alfa+\gamma_i\alfa)\xi_3
+\gamma_a\alfa\xi_6 +\gamma_i\alfa\xi_7+\delta_i\alfa(\xi_8-\xi_3)
+\frac{\beta\alfa(m-1)(\xi_1-\xi_2) S}{N},\\[1.2mm]
\Dl \xi_4(t') = k_2-(\gamma_a\alfa+\gamma_i\alfa)\xi_4
+\gamma_a\alfa\xi_6 +\gamma_i\alfa\xi_7+\delta_p\alfa(\xi_8-\xi_4)
+\frac{{\beta'}\alfa(m-1)(\xi_1-\xi_2) S}{N},\\
\Dl \xi_5(t') = \Dl \xi_7(t') = \Dl \xi_8(t') = 0,\\
\Dl \xi_6(t') =\gamma_r\alfa(\xi_7-\xi_6)+\delta_h\alfa(\xi_8-\xi_6)
-\frac{l\beta\alfa(m-1)(\xi_1-\xi_2)S}{N},
\end{cases}
\end{equation}
with $t'=t_f-t$. In turn, the minimality conditions of PMP establish
that the optimal controls $v^*$ and $m^*$ are given by
\begin{equation}
\label{eq:ext:cont}
\begin{split}
v^*(t) & =\min\left\{\max\left\{0,
\frac{(\xi_1(t)-\xi_7(t))S(t)}{2 k_3}  \right\},{v}_{\max}\right\},\\[2mm]
m^*(t) & =\min\left\{\max\left\{0,\frac{(\beta\alfa(I(t)+l H(t))
+{\beta'}\alfa P(t))(\xi_2(t)-\xi_1(t))S(t)}{2k_4 N}  \right\},{m}_{\max}\right\}.
\end{split}
\end{equation}

Finally, according to PMP, the following transversality conditions also hold true:
\begin{equation}
\label{eq:transversality}
\xi_i(t_f)=0, \quad i=1,\ldots,8.
\end{equation}

In Section~\ref{subsec:numres}, we use the necessary optimality conditions
\eqref{eq:co_states_fr2}--\eqref{eq:transversality} to numerically solve
the optimal control problem \eqref{cost-functional}--\eqref{omega:set},
both in classical and fractional cases.


\subsection{Numerical Results and Cost-Effectiveness of the Fractional Optimal Control Problem}
\label{subsec:numres}

Pontryagin Minimum Principle (PMP) is utilized to numerically solve the optimal control
problem as discussed in Section~\ref{subsec:control}. For that we use the
predict-evaluate-correct-evaluate (PECE) method of Adams--Basforth--Moulton
\cite{diethelm2005algorithms}, implemented by us in MATLAB. Firstly, we solve
System \eqref{sys:state} by the PECE procedure, with the initial values
for the state variables \eqref{ocp:ic} given, and a guess for
the controls over the time interval $[0, t_f ]$, thus obtaining
the values of the state variables. Analogously to \cite{Rosa},
a change in variables is employed in the adjoint system and
in the transversality conditions, obtaining a fractional initial
value problem (IVP) from \eqref{eq:co_states_fr2}--\eqref{eq:transversality}.
The resulting IVP is also solved with the PECE algorithm, and the values of the co-state
variables $\xi_i$, $i = 1, \ldots, 8$, are obtained. The controls are then updated
by a convex combination of the controls of  the previous iteration and the current values,
computed according to \eqref{eq:ext:cont}. This procedure is iteratively repeated
until the values of all the variables and the values of the controls are almost
coincident with the ones of the previous iteration, that is, until convergence is achieved.
The solutions of the classical model (i.e., the case $\alpha = 1$)
were successfully confirmed by a classical forward-backward scheme,
also implemented by us in MATLAB.

As estimated in \cite{cave2020covid}, we consider
that 10\% of infected individuals are super-spreaders.

According with \cite{lemos2020new}, the initial number
of asymptomatic persons is estimated to be
$$
\frac{I_0+P_0}{0.15}.
$$

In what follows, we assume that the total population is equal to 10,280,000 ($N$).
Based on real data obtained from \cite{DGS21} and on the  above assumptions,
the initial conditions are
\mbox{$R_0 = 278776$,} $E_0=92069$, $P_0= 68208\times 0.1$,
$I_0 = 68208\times 0.9$,
$A_0 = 68208/0.15$, $F_0 = 34$, $H_0 = 2366$ and
$S_0 = N - R_0 - E_0 - P_0 - I_0 - A_0 - H_0$.

The maximum number of effective daily vaccinated is estimated to be
30,000 (60,000 vaccine doses for two-dose vaccines), which corresponds
to 0.003 in percentage, and this is the value of $v_{\max}$.

In addition, we consider that the maximum rate of reduction in contacts
is the maximum value of function $m(t)$, $m_{\max}$, exhibited
in Figure~\ref{fig:function_m} and considered during  the modelling phase.

\begin{figure}[H]
\hspace{-1.5em}\includegraphics[scale=0.44]{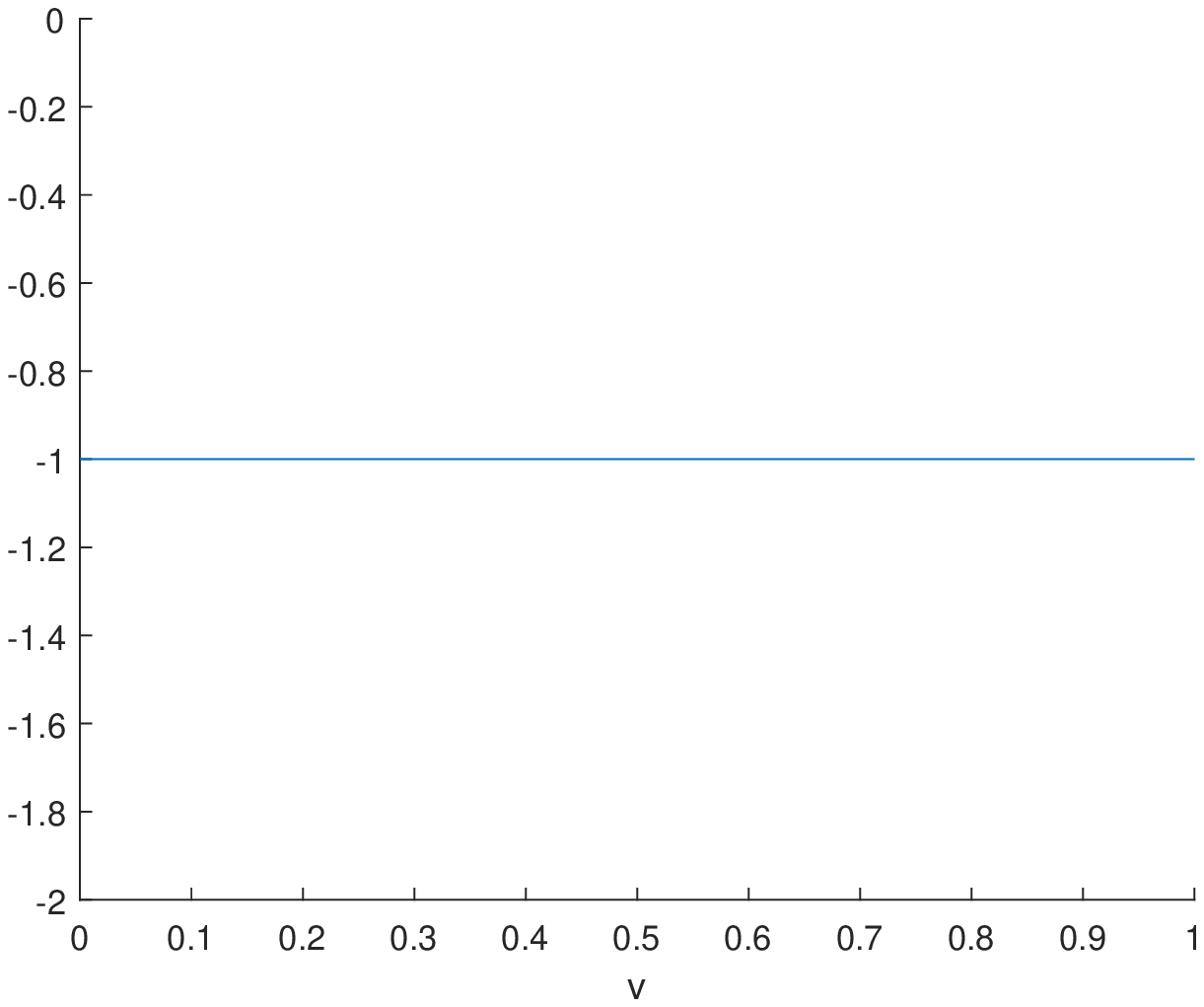}
\includegraphics[scale=0.44]{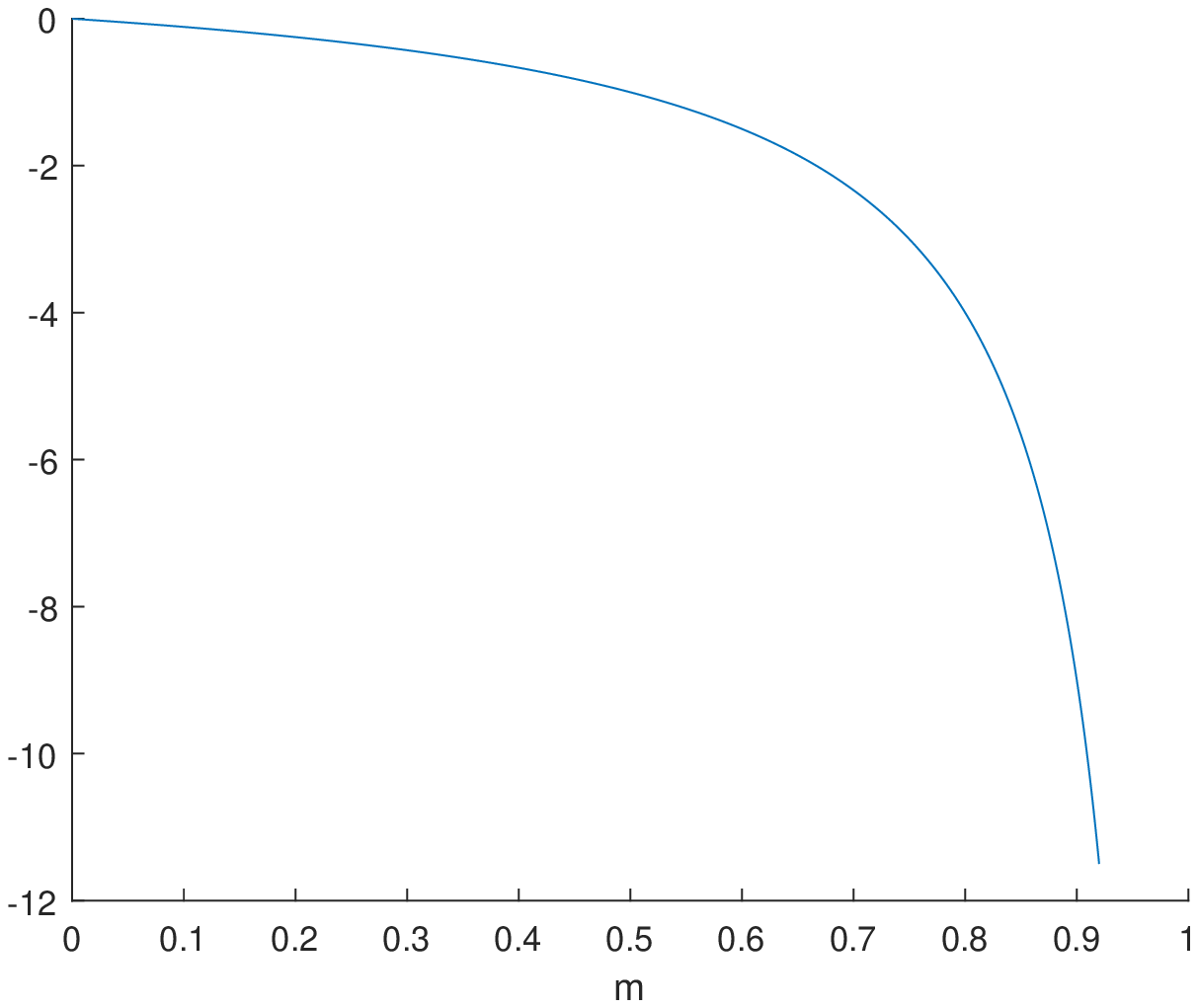}
\caption{Sensitivity index of the basic reproduction number
\eqref{r0_model_controls} with respect to the control
variables $v$ (left) and $m$ (right).}
\label{fig:Effic_controls}
\end{figure}

Following \cite{2019BpRL...14...27P}, the relevance of the two
control measures considered in the control of the disease
is calculated using the sensitivity index, as presented
in Definition~\ref{def:sentInd}. With this purpose,
the basic reproduction number of the model, now incorporating
the two controls, needs to be determined. For that, we look
to the controls as parameters. Using the next-generation matrix
method of \cite{van:den:Driessche}, we obtain that the
basic reproduction number is now given by
\begin{equation}
\label{r0_model_controls}
R_0 = \frac{(a_h + \gamma_a\alfa l)\beta\alfa(m-1)\rho_1}{a_h a_i v}
+ \frac{(a_h{\beta'}\alfa + \beta\alfa \gamma_a\alfa l)(m-1)\rho_2}{a_h a_p v},
\end{equation}
where $a_i= \gamma_a\alfa+\gamma_i\alfa+\delta_i\alfa$,
$a_p= \gamma_a\alfa+\gamma_i\alfa+\delta_p\alfa$,
and $a_h= \gamma_r\alfa+\delta_h\alfa$.
The sensitivity indices are presented as functions
of the control parameters in Figure~\ref{fig:Effic_controls},
using the parametric values from Table~\ref{tab:param}
and the classical derivative.
The plots of Figure~\ref{fig:Effic_controls} show the following:
(i) the curve of vaccination is constant and equal to $-1$,
meaning that  the vaccination program has a substantial impact, even for small
rates of application; (ii)~the curve of preventive measures
rapidly moves away from zero, meaning that a preventive programme
has a substantial impact only if its rate of application is high. Thus,
the use of masks and the limitation of individuals in confined spaces
(shops, schools, and other public spaces) have to be mandatory
and should be monitored by the authorities, when it is possible.

The weights of the cost functional \eqref{cost-functional} balance
the relative importance of quadratic control terms. Since
super-spreaders have a greater impact in the dynamics of the disease,
we consider that super-spreaders are more expensive to control
than regular infected individuals. Because the preventive measures
need a high rate of application to be effective, we consider that
preventive measures are more expensive than vaccination.
In our numerical experiments, weights are
$k_1 = 1$, $k_2 = 5$, $k_3 = 1$, and $k_4 = 10$.

\begin{figure}[H]
\hspace{-0.6em}\includegraphics[scale=0.44]{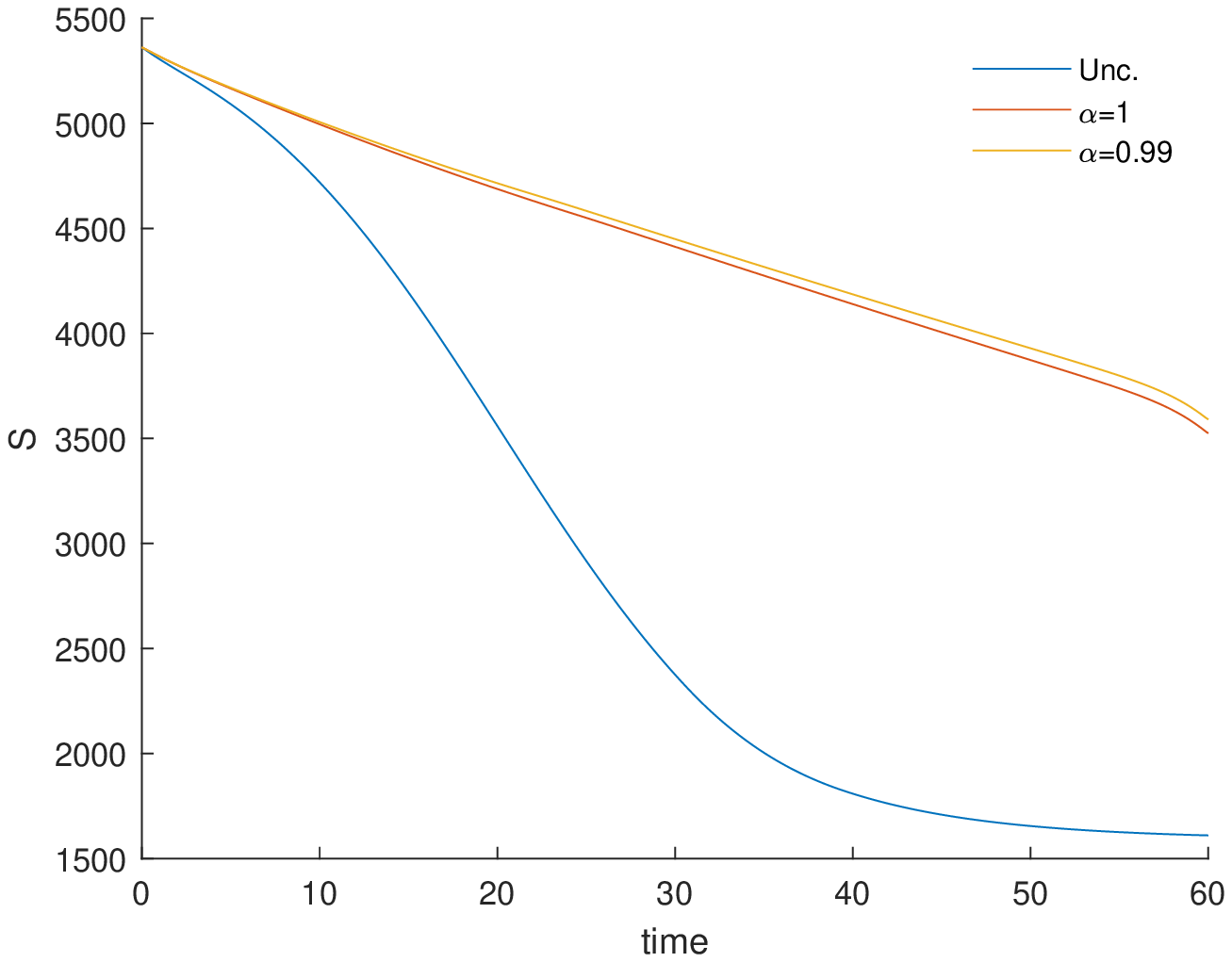}
\includegraphics[scale=0.44]{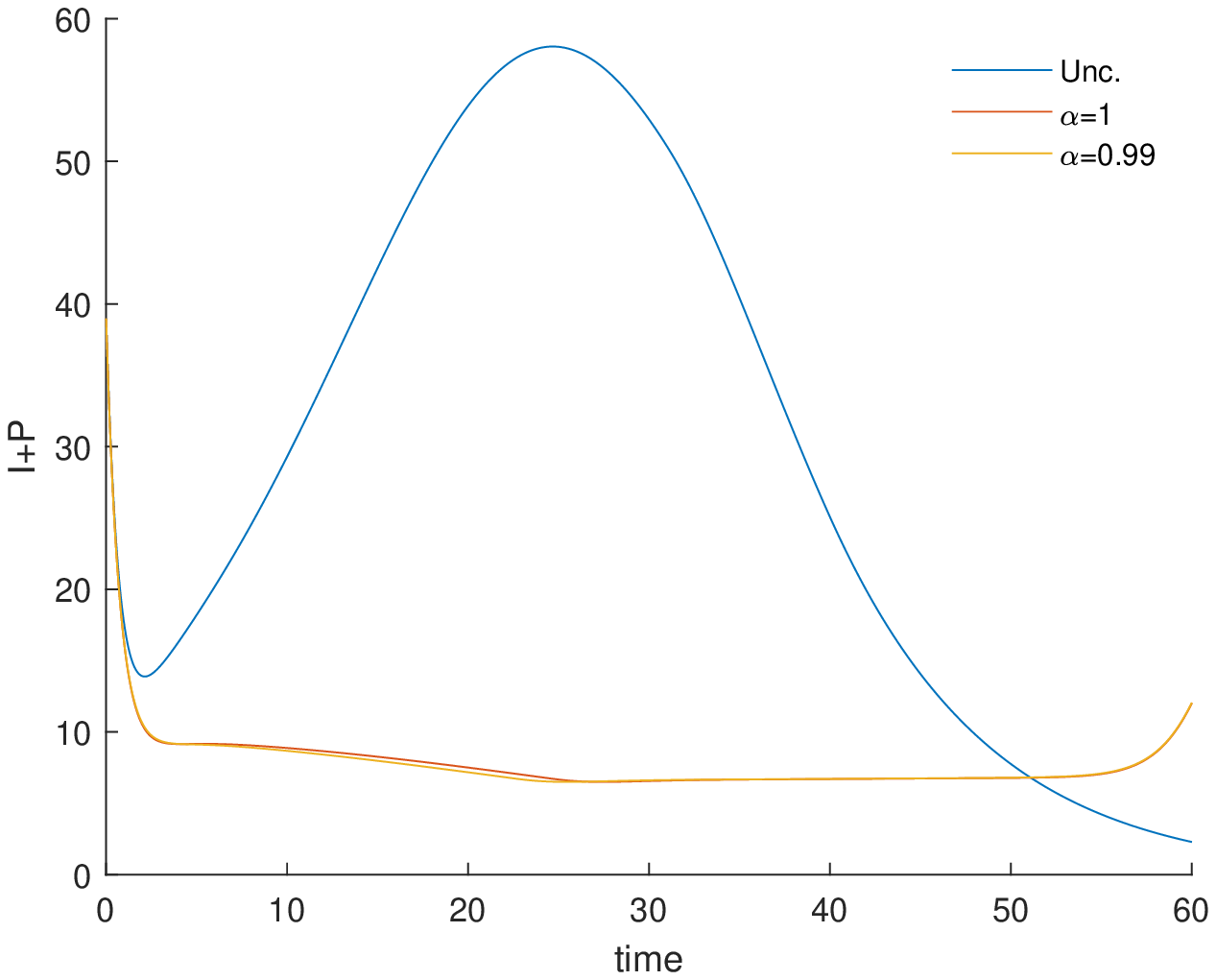}

\hspace{-0.6em}\includegraphics[scale=0.44]{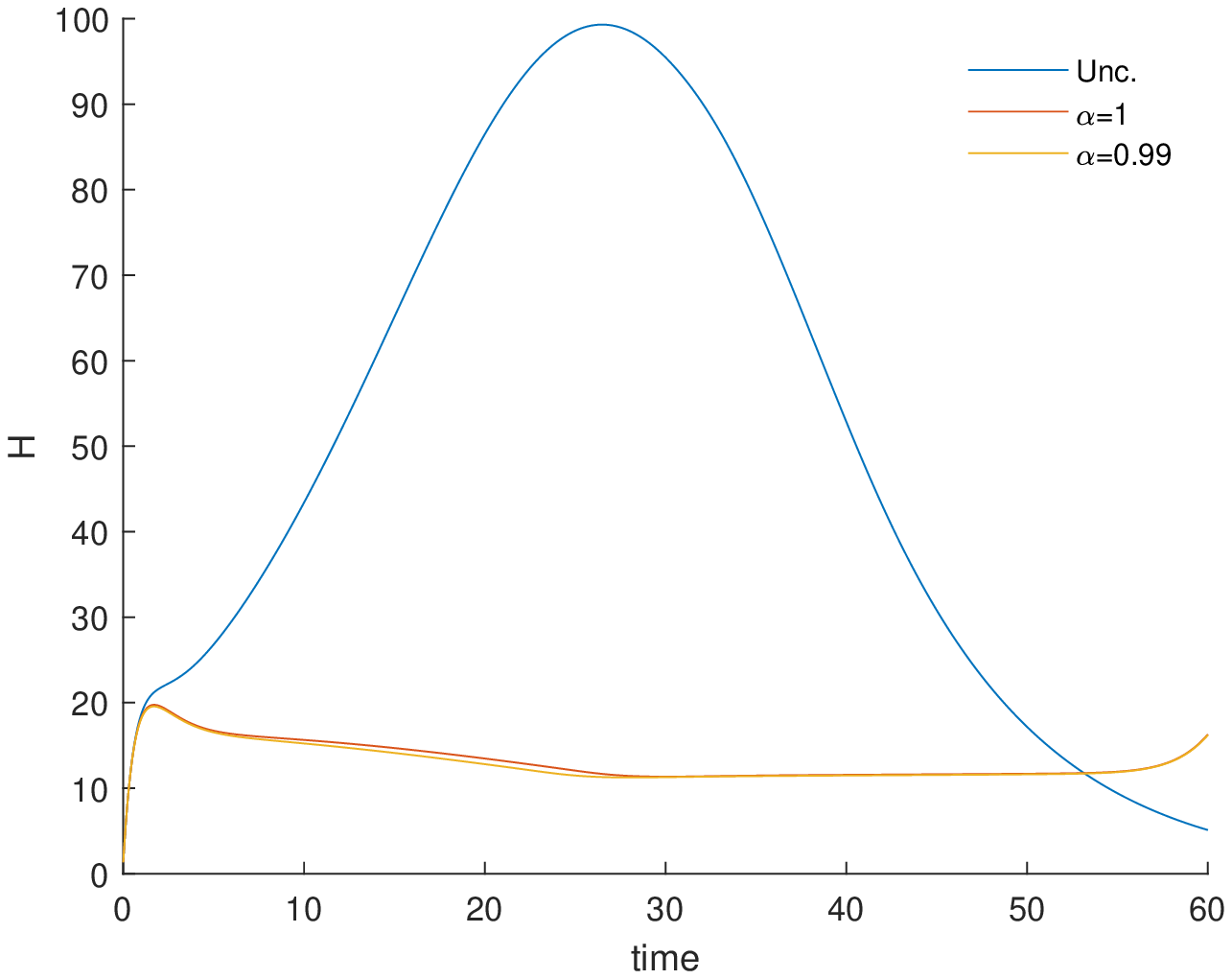}
~~~~~~~~\includegraphics[scale=0.44]{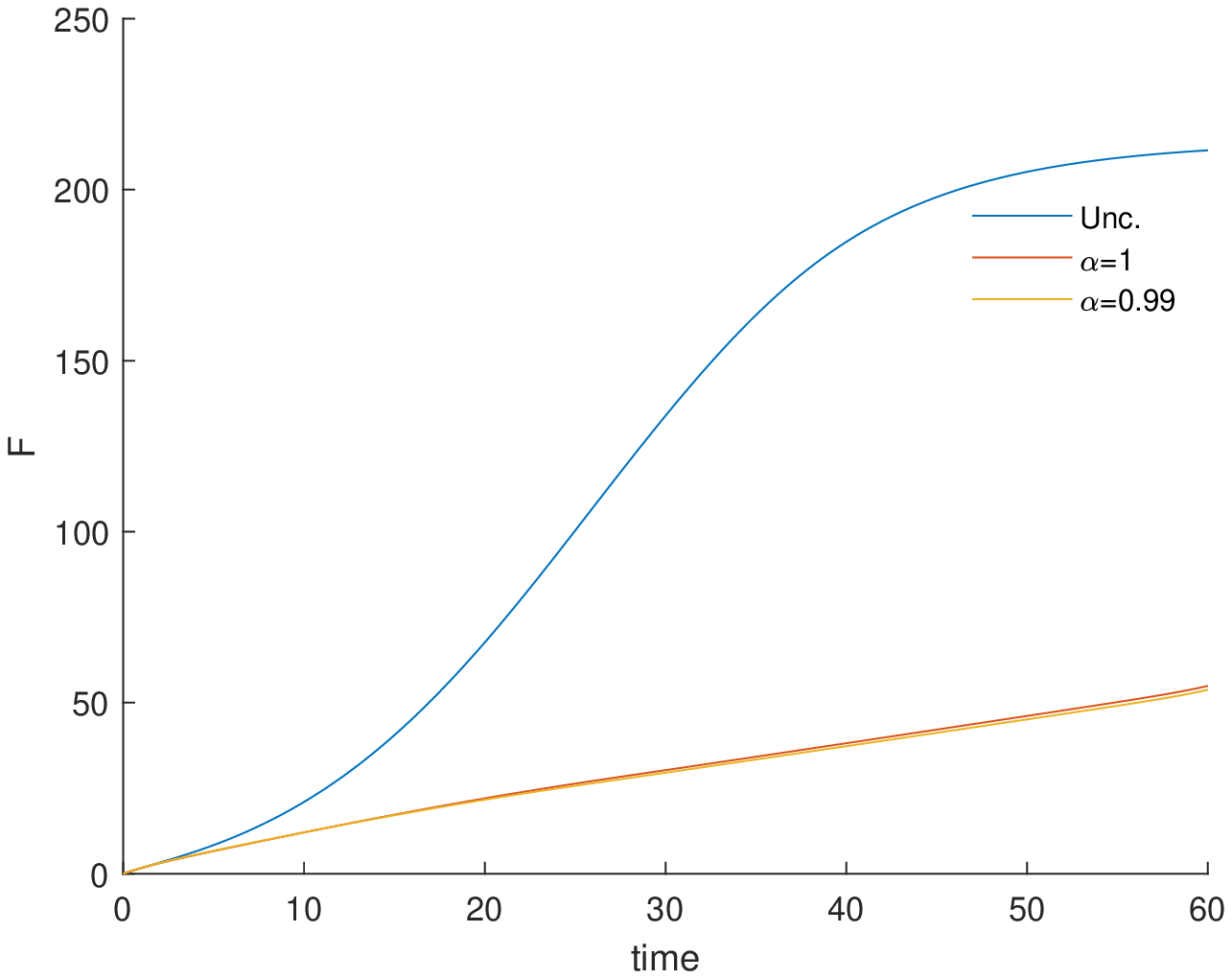}
\caption{Evolution of susceptible individuals $S$ (top left),
symptomatic infected individuals $I+P$ (top right),
hospitalized individuals $H$ (bottom left),
and fatalities $F$ (bottom right)
for the solutions of the uncontrolled model \eqref{Covid_model:exp}
and the optimal solutions of the FOCP \eqref{cost-functional}--\eqref{omega:set}
with fractional-order derivatives $\alpha=1.0$ and $\alpha=0.99$
and the parameter values of Table~\ref{tab:param}.}
\label{fig:S_E_var:alphas}
\end{figure}

In Figure~\ref{fig:S_E_var:alphas},
the trajectories of the fractional optimal control problem (FOCP)
for $\alpha=0.99$ are exhibited along with the solution
of the classical optimal control problem (i.e., with $\alpha=1$)
and the original (uncontrolled) model \eqref{Covid_model:exp}.

\begin{figure}[H]
\hspace{-0.7em}\includegraphics[scale=0.43]{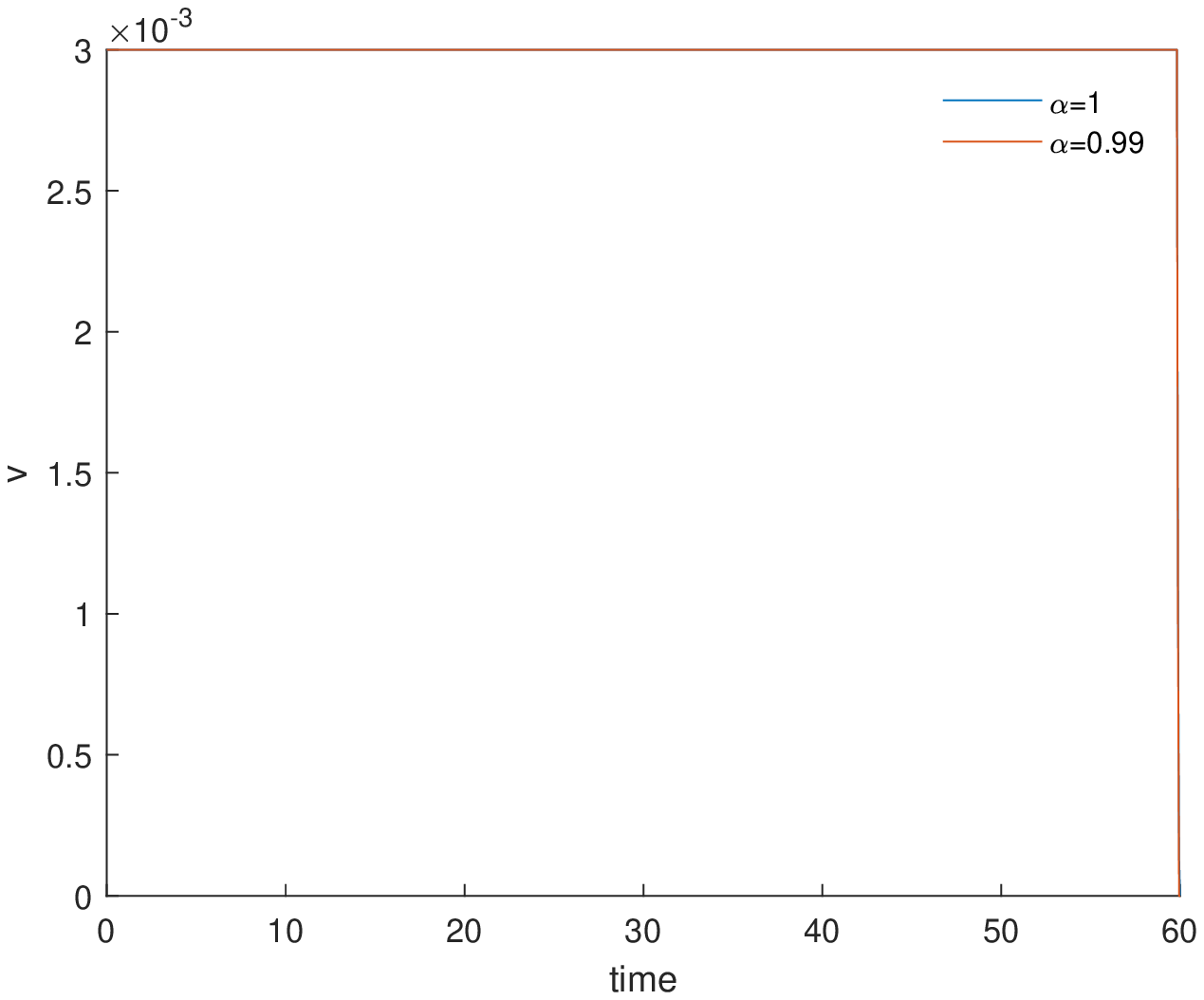}
\includegraphics[scale=0.43]{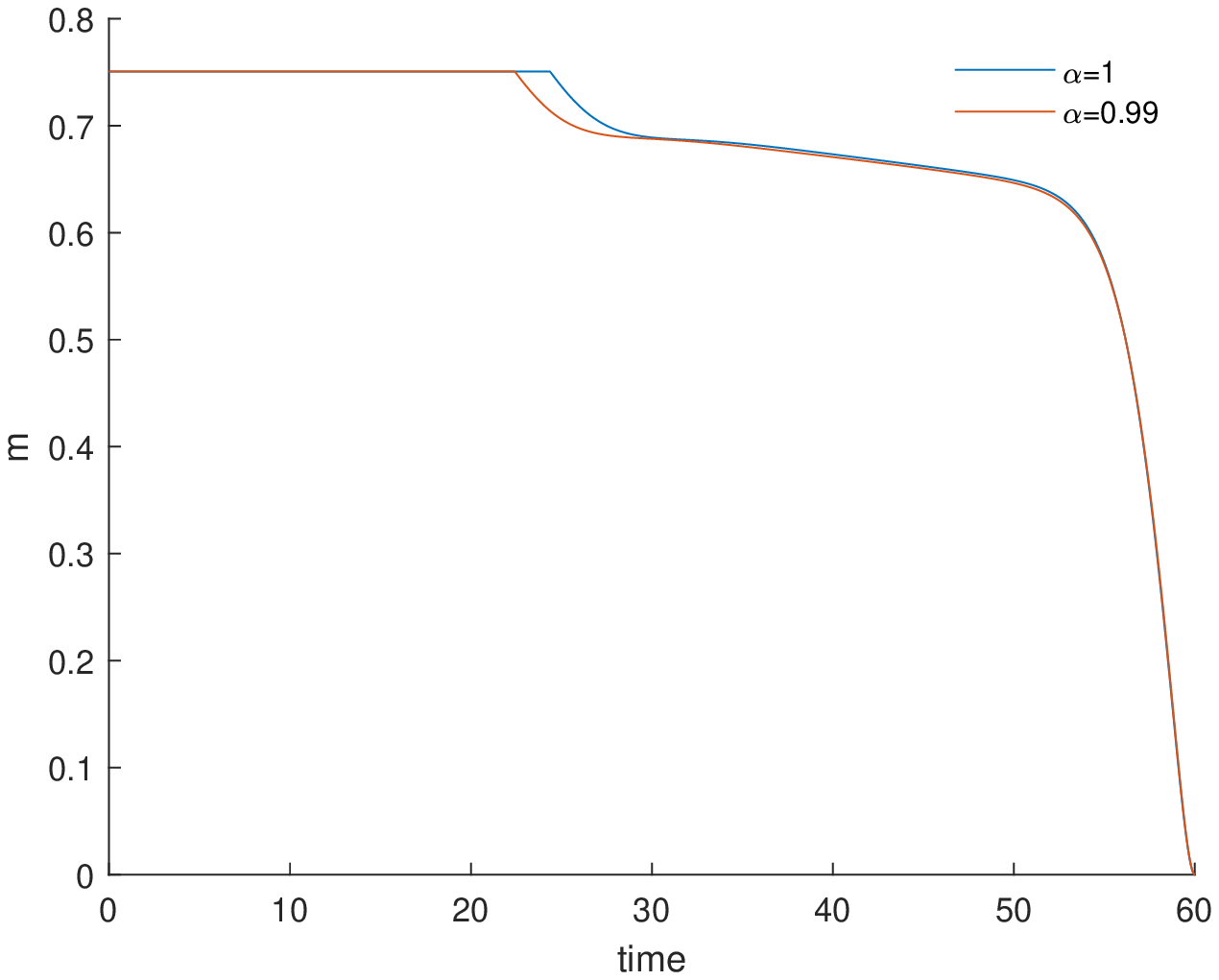}
\caption{Pontryagin controls $v$ (left) and $m$ (right)
for the FOCP \eqref{cost-functional}--\eqref{omega:set}
using the values in Table~\ref{tab:param} and fractional
order derivatives $\alpha=1.0$ and $\alpha=0.99$. The extremal
controls $v$ take their maximum value $v_{\max}$ almost everywhere.}
\label{fig:u_v_var:alphas}
\end{figure}

The corresponding Pontryagin controls are shown
in Figure~\ref{fig:u_v_var:alphas}.

We can see in Figures~\ref{fig:S_E_var:alphas} and \ref{fig:u_v_var:alphas}
that a change in the value of $\alpha$ corresponds to
variations in the state and control variables.
Moreover, comparing the solution of the original/uncontrolled model
with the solution of the optimal control problem obtained from the application
of the Pontryagin principle, we conclude that the considered control measures
are effective in the management of COVID-19.

\begin{figure}[H]
\hspace{-0.5em}\includegraphics[scale=0.6]{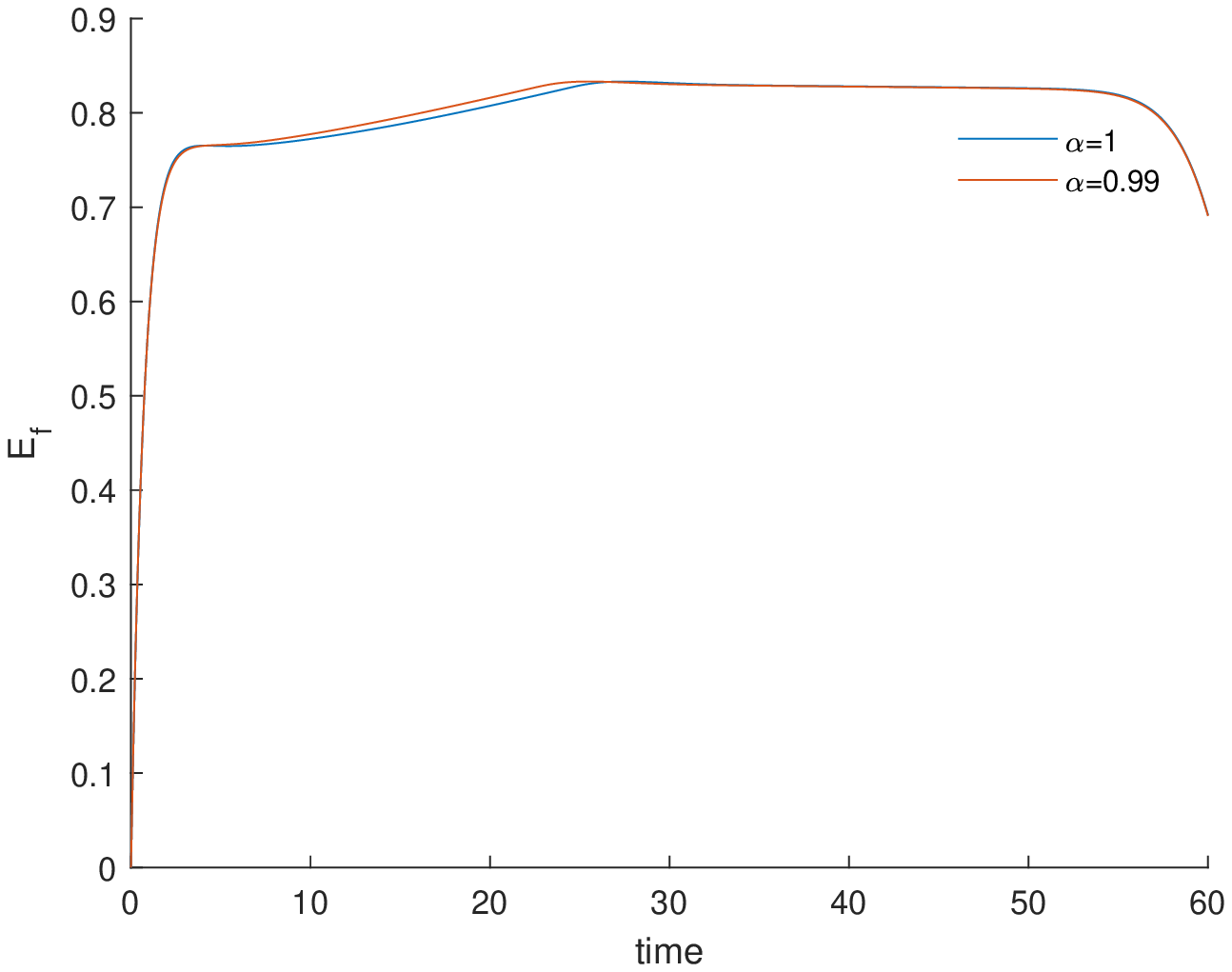}
\caption{Evolution of the efficacy function \eqref{efficacy_function}
for the FOCP \eqref{cost-functional}--\eqref{omega:set}
with values in Table~\ref{tab:param} and fractional-order
derivatives $\alpha=1.0$ and $\alpha=0.99$.}
\label{fig:Eff_varalpha}
\end{figure}

Figure~\ref{fig:Eff_varalpha} exhibits the efficacy function,
defined in \cite{rodrigues2014cost} by
\begin{equation}
\label{efficacy_function}
E_f(t)=\frac{i(0)-i^*(t)}{i(0)}
= 1-\frac{i^*(t)}{i(0)},
\end{equation}
where $i^*(t)=I^*(t)+P^*(t)$ is the optimal solution associated
with the fractional optimal control, and $i(0)=I(0)+P(0)$
is the correspondent initial condition.

The efficacy function $E_f(t)$ measures the proportional variation
in the number of infected individuals after the application
of the control measures, $\{v^*,m^*\}$, by comparing the number
of infectious individuals at time $t$ with its initial value $i(0)$.

To assess the cost and the effectiveness of the proposed fractional
control measures during the intervention period, some summary measures
are now presented.

The total cases averted by the intervention during
the time period $t_f$ is defined by
\begin{equation}
\label{eq:A}
AV=t_f i(0)-\int_0^{t_f}i^*(t)~dt,
\end{equation}
where $i^*(t)$ is the optimal solution associated
with the fractional optimal controls, and $i(0)$
is the correspondent initial condition \cite{rodrigues2014cost}. Note
that this initial condition is \textls[25]{obtained as the equilibrium proportion
$\overline{i}$ of System  \eqref{Covid_model:exp},
which is independent of time, so that
\mbox{$t_f i(0)=\int_0^{t_f}\overline{i}~dt$} represents the total infectious
cases over a given period of $t_f$ days.}

Effectiveness is defined as the proportion of cases averted
to the total cases possible under no intervention \cite{rodrigues2014cost}:
\begin{equation}
\label{eq:F}
\overline{F}=\frac{AV}{i(0) t_f}
=1-\frac{\displaystyle \int_0^{t_f}i^*(t)~dt}{i(0) t_f}.
\end{equation}

The total cost associated with the intervention
is defined in \cite{rodrigues2014cost} by
\begin{equation}
\label{eq:TCI}
TC=\int_0^{t_f} \left( C_1\,v^*(t)s^*(t)+C_2\,m^*(t)i^*(t)\right)~dt,
\end{equation}
where $s^*(t)=S^*(t)$, and $C_i$ corresponds to the per person unit cost
of the two possible interventions: (i) vaccination at time $t$
of susceptible individuals ($C_1$) and  (ii)   the  implementation of
preventive measures, such as the use of masks and  the physical distancing
of susceptible individuals~($C_2$).

Finally, the average cost-effectiveness ratio is given by
\begin{equation}
\label{eq:ACER}
ACER=\frac{TC}{AV}
\end{equation}
(see \cite{okosun2013optimal,rodrigues2014cost}).

Table~\ref{tab:efficacy} summarizes the presented
cost-effectiveness measures
\eqref{efficacy_function}--\eqref{eq:ACER}.
The results clearly show the effectiveness of the controls
in the reduction of COVID-19 infections and the advantage of
using the fractional model.
\begin{table}[H]
\caption{Summary of cost-effectiveness measures \eqref{efficacy_function}--\eqref{eq:ACER}
for classical ($\alpha =1$) and fractional ($0<\alpha <1$)
COVID-19 disease optimal control problems. Parameters according
to Table~\ref{tab:param} and $C_1=C_2=1$ in \eqref{eq:TCI}.}
\label{tab:efficacy}
\setlength{\tabcolsep}{8.5mm}{
\begin{tabular}{ccccc}
\toprule
\boldmath{$\alpha$} & \boldmath{$A$}  & \boldmath{$TC$}
& \boldmath{$ACER$} & \boldmath{$\overline{F}$}  \\[1mm] \midrule
0.99 & 1870.08 & 1116.43 & 0.596998 & 0.79967\\
1.0 & 1865.95 & 1114.53 & 0.597296 & 0.79791\\ \bottomrule
\end{tabular}}
\end{table}

\begin{figure}[H]
\includegraphics[scale=0.45]{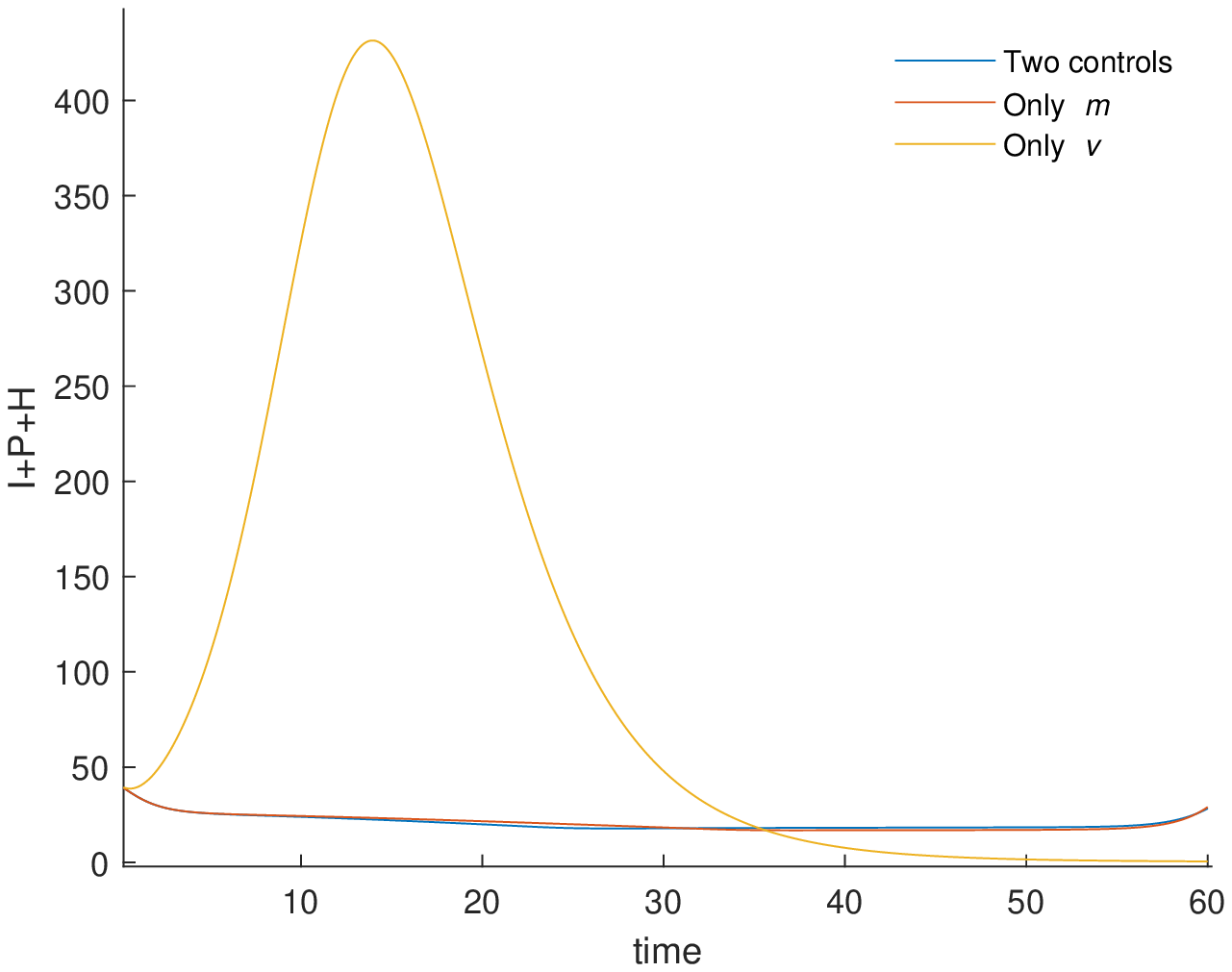}
\includegraphics[scale=0.45]{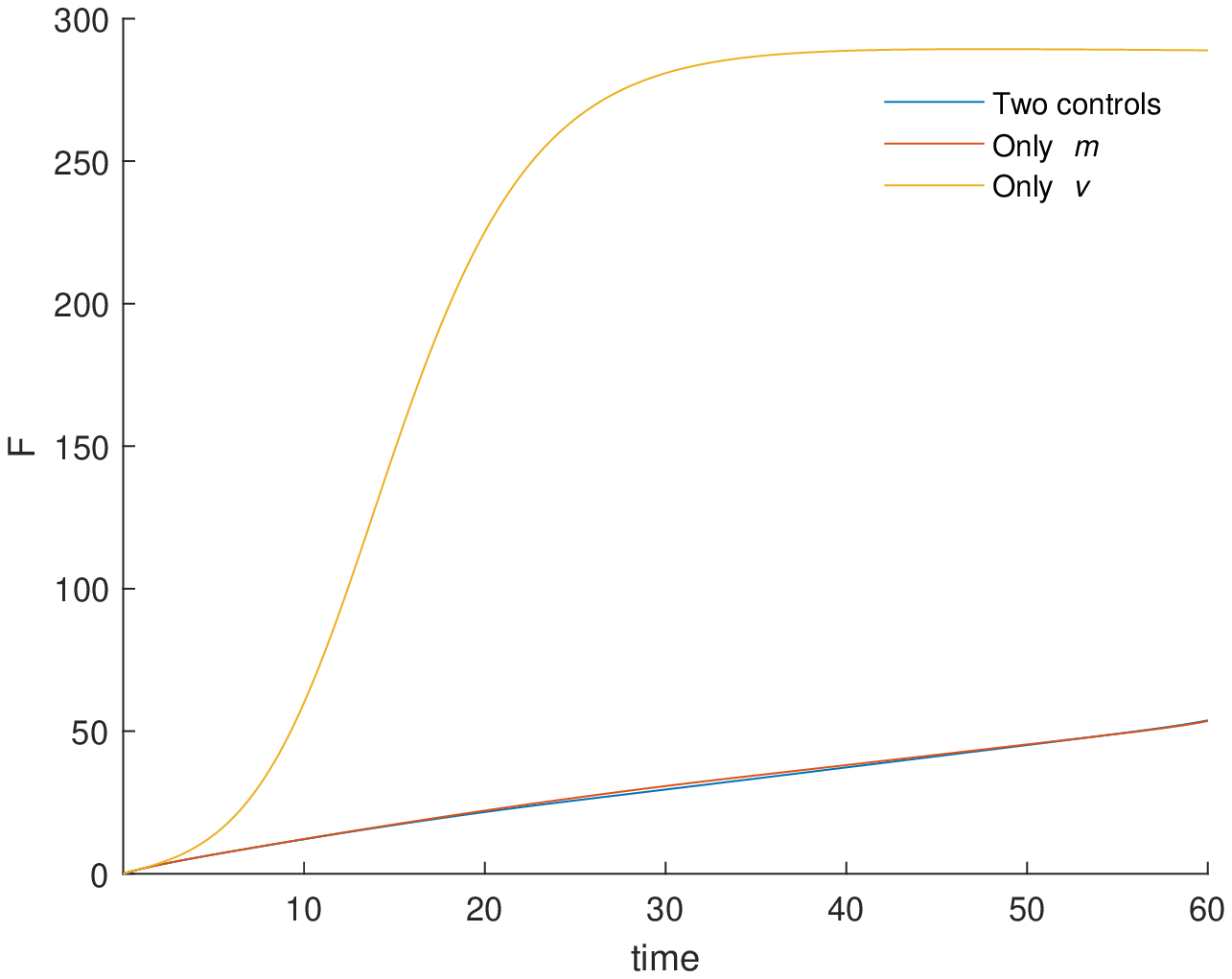}
\caption{Comparison of the solution of the FOCP
\eqref{cost-functional}--\eqref{omega:set} with the
fractional-order derivative $\alpha=0.99$, considering the two controls
with the two other cases where there is only one control used.
(left) Variation of infected individuals $I+P+H$.
(right) Evolution of fatalities $F$.}
\label{fig:alpha099:2controls_vs_1control}
\end{figure}

To conclude, Figure~\ref{fig:alpha099:2controls_vs_1control}
exhibits the dynamics of  the infected population, $I+P+H$, and of fatalities, $F$,
in three scenarios: (i)~when the \emph{two controls} are used;
(ii)~when \emph{only control $m$} is used ($v=0$); and
(iii)~when \emph{only control $v$} is used ($m=0$). Due to low vaccination rates,
considering $v=0$, one obtains almost the same solution as the one obtained using
both controls. On the other hand, erasing preventive measures leads to a serious healthcare problem.
Preventive measures are in this case more effective than vaccination.


\section{Conclusions}
\label{sec:conc}

A classical compartmental model with super-spreaders was firstly
proposed  and applied to provide an estimation of infected individuals
and deaths in Wuhan, China, in~\mbox{\cite{MR4093642}},
and this model was later extended to the fractional-order case in order
to include memory effects and better describe
the realities of Spain and Portugal
\cite{NDAIROU2021110652}. Unfortunately, the fractional-order
model of \cite{NDAIROU2021110652} is inconsistent, in the sense
that it does not satisfy appropriate time dimensions.
Here, the fractional model of \cite{NDAIROU2021110652} was
corrected and then used, for the first time in the literature,
to model the third wave of  the COVID-19 pandemic in Portugal,
which occurred between 27 December 2020 and 16 February 2021.
Our data fitting consisted in the minimization of the $l_2$ norm
of the difference between the real values reported from the Health Authorities
and  the predictive cases of COVID-19 infection given by our model \eqref{Covid_model:exp},
showing that, in terms of absolute errors, the fractional-order model
is better than the classical integer-order one. Another advantage
of the fractional-order model was found from a sensitivity analysis,
measuring the importance of each parameter in   COVID-19 transmission,
which allowed us to show that the sensitivity of the model decreases
in absolute value with the decrease in the fractional-order $\alpha$
of differentiation, i.e., the fractional-order model is less sensitive
to disturbances in the parameters than the classical one.
Moreover, we introduced into the corrected model the use of vaccination
and preventive control measures, investigating the use
of fractional-order optimal control theory to
minimize the number of COVID-19-infected
individuals and reduce the associated costs.
A post-optimal cost-effectiveness analysis
has shown the effectiveness of  controls
to combat COVID-19 and the advantage
of using the fractional-order model.
Finally, it is shown that preventive measures
are essential in the control of the pandemic.

The model investigated here is deterministic.
For future work, it would be interesting to take into account
the effect of noise. For integer-order models,
one can refer to the works~\mbox{\cite{Zine1,Zine2}},
where some COVID-19 stochastic differential equations
are proposed that take into account noises derived
from environmental fluctuations. Fractional stochastic models
for COVID-19 are scarce and still need further investigations.

\vspace{6pt}
\authorcontributions{Conceptualization, S.R. and D.F.M.T.; software, S.R.;
validation, D.F.M.T.; formal analysis, S.R. and D.F.M.T.; investigation,
S.R. and D.F.M.T.; writing---original draft preparation, S.R. and D.F.M.T.;
writing---review and editing, S.R. and D.F.M.T. All authors have read
and agreed to the published version of the manuscript.}

\funding{This research was funded by The Portuguese Foundation for Science and Technology
(FCT---Funda\c{c}\~{a}o para a Ci\^{e}ncia e a Tecnologia) through
projects UIDB/50008/2020 (Rosa) and UIDB/04106/2020 (Torres).}

\institutionalreview{Not applicable.}

\informedconsent{Not applicable.}

\dataavailability{The data supporting  the reported results
are public and can be found in~\cite{DGS21,Data}.}

\acknowledgments{The authors are very grateful
to two referees for several constructive remarks
and questions that helped them to improve the manuscript.}

\conflictsofinterest{The authors declare that there is no conflict of interest.
The funders had no role in the design of the study, in the collection,
analyses, or interpretation of data, in the writing of the manuscript,
or in the decision to publish the~results.}

\end{paracol}


\reftitle{References}


\end{document}